\documentclass{amsproc}
\UseRawInputEncoding
\usepackage{graphicx}

\usepackage{enumerate,hyperref} 
\usepackage[version=3]{mhchem}
\usepackage[all]{xy}
\usepackage{adjustbox}
\usepackage{arydshln}
\usepackage{marginnote}
\usepackage{multirow}
\usepackage{resizegather}
\usepackage[dvips,pdftex]{color}
\usepackage{colortbl}
\usepackage{pgf,tikz}
\usepackage{url}
\usepackage{longtable}
\usepackage[normalem]{ulem}

\newtheorem{theorem}{Theorem}[section]

\newtheorem{proposition}[theorem]{Proposition}

\newtheorem{conjecture}[theorem]{Conjecture}

\theoremstyle{definition}

\newtheorem{example}[theorem]{Example}

\theoremstyle{remark}

\numberwithin{equation}{section}

\newcommand{\R}{\mathbb{R}}
\newcommand{\C}{\mathbb{C}}

\newcommand{\Z}{{\mathbb Z}}
\newcommand{\Q}{{\mathbb Q}}

\newcommand{\mC}{\mathcal{C}}
\newcommand{\mG}{\mathcal{G}}
\newcommand{\mR}{\mathcal{R}}
\renewcommand{\k}{\kappa}

 \DeclareMathOperator{\diag}{diag}
 \DeclareMathOperator{\im}{Im} 

\newcommand{\New}{{\rm New}}

\newcommand{\invtPoly}{\mathcal{P}}

\newcommand{\idealB}{\mathfrak{B}}
\newcommand{\idealI}{\mathfrak{T}}  

\newcommand{\idealP}{\mathfrak{P}}

\begin{document}

\title[Chemical reaction networks]{From chemical reaction networks to algebraic and \\ polyhedral geometry -- and back again}


\author{Elisenda Feliu}
\address{Department of Mathematical Sciences, University of Copenhagen, Universitetsparken 5, 2100 Copenhagen, Denmark} 
\email{efeliu@math.ku.dk}
\thanks{The first author was supported by the Independent Research Fund of Denmark.}

\author{Anne Shiu}
\address{Department of Mathematics, Texas A\&M University.  TAMU Mailstop 3368, College Station TX 77843-3368, USA}
\email{annejls@tamu.edu}
\thanks{The second author was supported by the NSF (DMS-1752672).}



\date{\today}



\begin{abstract}
This is a chapter for a book in honor of Bernd Sturmfels and his contributions.  
We describe the contributions by Bernd Sturmfels and his collaborators in
harnessing algebraic and combinatorial methods for analyzing chemical
 reaction networks. Topics explored include the steady-state variety,
   counting steady states, and the global attractor conjecture.  
We also recount some personal stories that highlight Sturmfels's long-lasting impact on this research area.
	\end{abstract}

\maketitle

	\section*{Introduction}

How did Bernd Sturmfels get into the topic of chemical reaction networks?  The story begins at the Mathematical Sciences Research Institute (MSRI) in 2003, during a yearlong program on `Commutative Algebra'\footnote{In fact, some seeds were planted earlier: 
In the mid-1990s, Sturmfels had some initial discussions with Karin Gatermann, and separately heard from his Computer Science colleague Alistair Sinclair about discrete-time versions of chemical reaction systems~\cite{rabinovich1992quadratic}.}.  Karin Gatermann  was visiting and presented new ideas about sparse polynomial systems arising from chemical reaction networks~\cite{Karin01, Karin02}.  One day, Gatermann and Sturmfels were talking at a blackboard, when Alicia Dickenstein   walked by.  Sturmfels invited her to join the conversation, and, as they say, the rest is history.


History, however, is sometimes bittersweet.  Within a few years of her visit to MSRI, Gatermann died of cancer.  
To honor her, friends and colleagues organized a special issue of {\em Journal of Symbolic Computation}.   
Sturmfels proposed to submit a research article on chemical reaction networks to this issue, in part to carry on Gatermann's legacy and research in this area.  This was co-authored with Dickenstein, Gheorghe Craciun, and a Ph.D.\ student of Sturmfels -- Anne Shiu, an author of this chapter~\cite{TDS}.

\begin{figure}[htbp]
\begin{center}
  \includegraphics[width=4.5in]{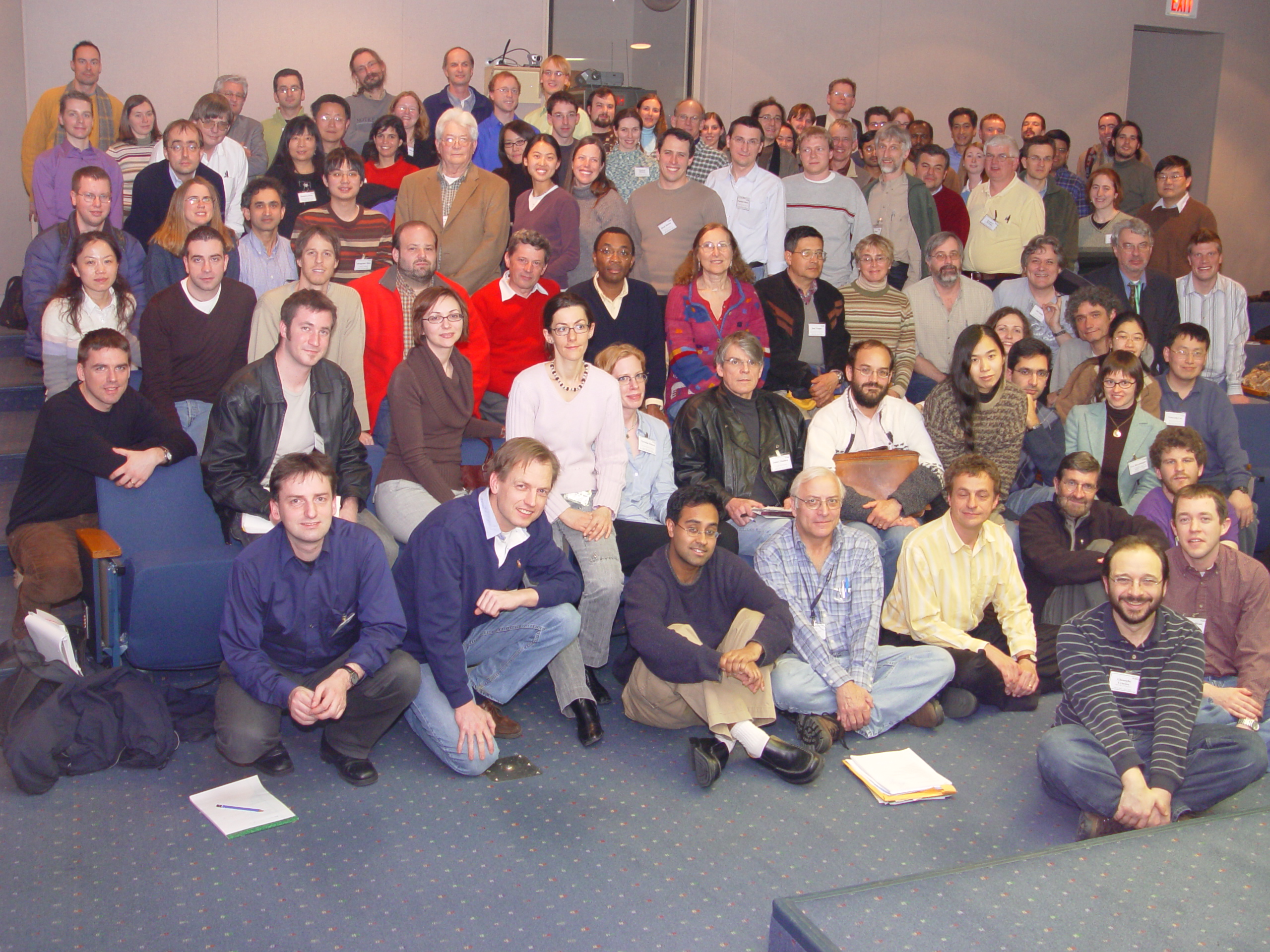}
\caption{Group photo at `Applications in Biology, Dynamics, and Statistics' workshop at IMA in March 2007.  Sturmfels is seated in the first row, fifth from the left.
Dickenstein is in the third row, seventh from the left; and Craciun is the unique person in the zeroth row.
}
\label{fig:IMA}
\end{center}
\end{figure}

Much of that first article on chemical reaction networks was developed at a 2007 Institute for Mathematics and its Applications (IMA) workshop on `Applications in Biology, Dynamics, and Statistics'; see the group photo in Figure~\ref{fig:IMA}.  The four authors met daily in a cafe at 7:00 am -- those who know Sturmfels know that he wakes up early!

In this chapter, we recount additional personal stories behind Sturmfels's research impact in chemical reaction networks.  Indeed, while Sturmfels is an author of only four articles in this area
\cite{TDS,ShiuSturmfels,CKLS,case-study}, a true measure of his impact is through the many people he brought together in collaborations.  This chapter also highlights the many ideas -- often borrowed from adjacent research areas -- that he was the first to bring into this field.  In fact, a theme running through Sturmfels's research on chemical reaction networks is summarized in the stated aim of~\cite{case-study}: ``to demonstrate how biology can lead to interesting questions in algebraic geometry and to apply state-of-the-art techniques from computational algebra to biology.'' 

This chapter is organized as follows.  In Section~\ref{sec:background}, we introduce chemical reaction networks and give an overview of the main questions in this research area.  Section~\ref{sec:complex-bal} concerns chemical reaction systems having a special type of steady state, namely, complex-balanced steady states, and their associated toric structures.  We also discuss the global attractor conjecture, including how Sturmfels gave the name to this conjecture.  
Siphons and their relation to boundary steady states are the focus of Section~\ref{sec:siphon}, 
followed by steady-state invariants and their relation to matroids in Section~\ref{sec:matroid}.
Next, Section~\ref{sec:mss} describes polyhedral methods for assessing whether a network admits multiple steady states, and then 
Section~\ref{sec:counting} explains how the mixed volume and techniques from numerical algebraic geometry are used to count or bound the number of steady states.  
Section~\ref{sec:attainable} concerns convex hulls of trajectories of chemical reaction systems, and we end with a discussion in Section~\ref{sec:end}.
	
	\section{Background} \label{sec:background}

	\subsection{Reaction networks and associated models}	
		In this section we  briefly introduce the mathematical setting of the theory of reaction networks.  This theory was initiated by chemical engineers who wanted to model the time-dependent evolution of chemical reactions.
		
		\medskip
		\noindent
	\paragraph{\bf Reaction networks. }	
	Reaction networks model the interactions among a finite set of objects, called {\em species}.  
	In this chapter, we let $\mathcal{S}=\{X_1,\dots,X_n\}$ denote the set of species.  
	A \emph{complex} is a finite linear combination of the elements in $\mathcal{S}$ with non-negative integer coefficients. 
We identify each complex with the vector it defines in $\Z^n$, 
so that a complex can be written in two forms:
	\[ y = \sum_{i=1}^n \alpha_i X_i \qquad \textrm{or}\qquad
	y=(\alpha_1,\dots,\alpha_n)\in \Z^n_{\geq 0}~.\]

	A \emph{reaction network} $\mathcal{G}=(\mathcal{C},\mathcal{R})$ is a digraph whose nodes are complexes. Edges are referred to as \emph{reactions}.  Throughout this chapter, the number of reactions is denoted by $r$ (so, $|\mathcal R| = r$) and the number of complexes by $m$ (so, $|\mathcal C| = m$). 
	
	\smallskip
	For example, the following  network is McKeithan's `kinetic proofreading' model of T-cell signal transduction~\cite{T-cell, Sontag01}: 
\begin{equation}\label{eq:mckeithan}
 \begin{xy}<8mm,0cm>:
(0,1.3)                  ="A+B"  *+!D{ X_1+X_2 }  *{};
(-2,0)                  ="D"  *+!R{X_4}  *{};
 (2,0)                  ="C"  *+!L{X_3}  *{};
    {\ar "A+B"+(-0.15,0)*{};"C"+(-0.3,.15)*{}  };
   {\ar "C"+(0,0.15)*{};"A+B"+(0.15,0)*{}};      
   {\ar "C"*{};"D"+(0.15,0)*{}};
  {\ar "D"+(0.15,0.15)*{};"A+B"-(0.15,0.15)*{}};\end{xy}
  \end{equation}
In this network, which we call the {\em McKeithan network}, $X_1$ represents the T-cell receptor, $X_2$ is the major histocompatibility complex (MHC) of an antigen-presenting cell, $X_3$ represents $X_1$ and $X_2$ bound together, and $X_4$ is the activated form of $X_3$.
  The set of complexes (nodes)  and reactions (edges) are
\begin{equation}\label{eq:mckeithanCR}
\begin{aligned}
\mC& =\{ X_1+X_2,~ X_3,~ X_4\}, \\ 
\mR & =\{ X_1+X_2\ce{->}X_3,~ X_3\ce{->} X_1+X_2,~X_3\ce{->} X_4, ~X_4 \ce{->} X_1+X_2\}. 
\end{aligned}
\end{equation}
Hence, there are $n=4$ species, $m=3$ complexes, and $r=4$ reactions.

		\medskip
As mentioned earlier, reaction networks represent interactions among species.  
These species may arise in chemistry, biochemistry, or even in population biology or epidemiology.
Therefore, species can be chemical compounds, proteins (as in the network in \eqref{eq:mckeithan}), animals, classes of people, and so on. 
For instance, the classical Lotka-Volterra model of prey-predator interaction has the following underlying reaction network: 
	\[ 0 \ce{->} P, \qquad P+U \ce{->} 2U, \qquad U \ce{->} 0,\]
where $P$ is the prey and $U$ the predator.  
As another example, the classical SIR model in epidemiology builds on the following reaction network:
	\[ S +I \ce{->} 2I, \qquad I \ce{->} R,\]
where $S$, $I$, $R$ refer to susceptible, infected, and removed individuals.

			\medskip
		\noindent
	\paragraph{\bf Dynamical systems. }	
Given a reaction network, one aims at modeling the time evolution of the species' abundances. 
This chapter focuses on deterministic (rather than stochastic) models with continuous (rather than discrete) time.  Hence, 
we consider the concentrations of the species, and we model their time evolution using systems of ordinary differential equations (ODEs). 

Formally, we let $x_i(t)$ denote the concentration of species $X_i$ at time $t$, and we consider an ODE system of the form
\begin{equation}\label{eq:ode}
\frac{dx_i}{dt} = \sum_{y\rightarrow y'\in \mR} v_{y\rightarrow y'}\big(x(t)\big) (y_i'-y_i), \qquad i=1,\dots,n.
\end{equation}
Here $x(t)=(x_1(t),\dots, x_n(t))$,  
and $(y_i'-y_i)$ represents the net production of species $X_i$ resulting from reaction $y \to y'$.   Also, 
$v_{y\rightarrow y'}$ is a continuously differentiable function in $x_1,\dots,x_n$. A choice of functions 
$v_{y\rightarrow y'}$ for all reactions $y\rightarrow y'$ is referred to as a choice of
\emph{kinetics}. 
A widespread choice of kinetics is \emph{mass-action kinetics}, in which 
 the rate of each reaction is proportional to the product of the concentrations of the reactant species of the reaction.  That is, each $v_{y\rightarrow y'}$ is a monomial of the following form:
\begin{equation}\label{eq:mass-action}
 v_{y\rightarrow y'}(x) =  \k_{y\rightarrow y'} \, x_1^{y_1} \cdots x_n^{y_n}, \qquad \textrm{where }\ \k_{y\rightarrow y'} >0.
\end{equation}
The proportionality constant $\k_{y \to y'}$ is called the \emph{reaction rate constant}. In this case, the right-hand side of each ODE in system~\eqref{eq:ode} is a polynomial, and this is where the world of applied algebra comes into play. 

The ODE system \eqref{eq:ode} admits a decomposition as the product of a matrix with the vector of kinetics, which we describe next. 
Choose an ordering of the reaction set $\mathcal R$, and construct the \emph{stoichiometric matrix} $N\in \R^{n\times r}$, in which the $i$-th column is the vector $y\rightarrow y'$ if the $i$-th reaction is $y\rightarrow y'$. 
Also, let 
$v(x)$ denote the length-$r$ vector in which the entry corresponding to the vector $y \to y'$ is $v_{y\rightarrow y'}(x)$.
Now, after omitting reference to $t$, we rewrite \eqref{eq:ode} in matrix-vector form:
\begin{equation}\label{eq:odeN}
\frac{dx}{dt}= N v(x).
\end{equation}

Next, we describe other 
useful decompositions of the ODE system \eqref{eq:ode}, for the case of mass-action kinetics.  Choose an ordering of the set of complexes: $\mC = \{y_1,\dots,y_m\}$. Let $Y\in \R^{n\times m}$ denote the matrix whose columns are the complexes, and let $x^Y= ( x^{y_1},\dots,x^{y_m})^{\top} 
\in \R^{m}$. 
 (Recall that  $x^{y_i}$ denotes the monomial $x_1^{y_{i1}} \cdots x_m^{y_{im}}$).  
Let $C_{\mathcal{G}}\in \R^{m\times r}$ denote the incidence matrix of $\mG$, and let $K_\k\in \R^{r\times m}$ denote the matrix whose $(i,j)$-th entry is $\k_{y\rightarrow y'}$ if $y\rightarrow y'$ is the $i$-th reaction and $y$ the $j$-th complex. Then $A_\k:=C_{\mathcal{G}} K_\k\in \R^{m\times m}$ is the negative of the Laplacian matrix of the labeled digraph obtained from $\mathcal{G}$ by labelling each reaction $y\rightarrow y'$ with its reaction rate constant $\k_{y\rightarrow y'}$. 

With this notation in place, and after observing that $N=Y C_{\mathcal{G}}$, the right-hand side of \eqref{eq:ode} with mass-action kinetics 
admits the following two decompositions:
\begin{equation}\label{eq:odeYCG}
 N v(x) ~=~ Y C_{\mathcal{G}} K_\k\,  x^Y ~=~ Y A_\k\,  x^Y. 
 \end{equation}

 \begin{example}\label{ex:mc}
Let us find the ODE system under mass-action kinetics for the McKeithan network \eqref{eq:mckeithan} using the objects just introduced. 
We choose the orderings of the complexes and reactions given in \eqref{eq:mckeithanCR}, and denote by  $\k_i$   the reaction rate constant of the $i$-th reaction. 
This gives rise to the following objects:  
{\small \begin{align*}
N  & = \left[\begin{array}{rrrr}
-1 & 1 & 0 & 1 \\
-1 & 1 & 0 & 1 \\
1 & -1 & -1 & 0 \\
0 & 0 & 1 & -1
\end{array}\right], &  
v(x)& = ( \k_1x_1x_2,~ \k_2 x_3,~ \k_3x_3,~ \k_4x_4)^{\top}, \\
Y &= \begin{bmatrix}
1 & 0 & 0 \\ 1 & 0 & 0 \\ 0 & 1 & 0 \\ 0 & 0 & 1
\end{bmatrix}, &
 C_\mG &= \left[\begin{array}{rrrr}  
-1 & 1 & 0 & 1 \\
1 & -1 & -1 & 0 \\
0 & 0 & 1 & -1
\end{array}\right], 
\\
 K_\k & =  \left[\begin{array}{ccc}  
\k_1  & 0 & 0 \\
0 &  \k_2 & 0 \\
0 &  \k_3 & 0 \\
0 & 0 &  \k_4
\end{array}\right], 
& 
A_\k &=  \left[\begin{array}{ccc}  
-\k_1 & \k_2 & \k_4  \\
\k_1 & -\k_2-\k_3 & 0  \\
0 & \k_3 & -\k_4  
\end{array} \right].
\end{align*}}%
Using \eqref{eq:odeN} or \eqref{eq:odeYCG}, the ODE system describing the time evolution of the concentrations $x_1,x_2,x_3,x_4$ is as follows:
\begin{equation}\label{eq:odeex}
 \begin{aligned}
 \tfrac{dx_1}{dt} &= - \k_1 x_1x_2 + \k_2 x_3 + \k_4 x_4 &   \quad \tfrac{dx_3}{dt} &=  \k_1 x_1x_2 - \k_2 x_3 - \k_3 x_3  \\
  \tfrac{dx_2}{dt} &= - \k_1 x_1x_2 + \k_2 x_3 + \k_4 x_4 &   \quad  \tfrac{dx_4}{dt} &= \k_3 x_3 - \k_4 x_4.
 \end{aligned}
 \end{equation}
 \end{example}	
 
 \medskip

 From the ODEs~\eqref{eq:odeex} of the McKeithan network, we easily see that 
 \begin{align} \label{eq:conservation-McK}
  \tfrac{dx_1}{dt}  +  \tfrac{dx_3}{dt}  +  \tfrac{dx_4}{dt}  =0,
  \end{align}
 and hence the sum $x_1(t)+x_3(t)+x_4(t)$ is a constant function in $t$ for every solution of the ODE system. Its value equals the sum at the initial condition $(x_1(0),x_2(0),x_3(0),x_4(0) )$. 

The conservation relation~\eqref{eq:conservation-McK} corresponds to the sum of the first, third, and fourth row of $N$ being zero. Indeed, every vector in the left kernel of $N$ gives such a relation. For the McKeithan network, 
one additional independent such relation exists, for instance, $x_2(t)+x_3(t)+x_4(t)$ is also constant. 
 
In general, linear conservation relations arise as follows.  
From~\eqref{eq:odeN}, we see that the derivatives $\frac{dx}{dt}$ belong to the vector subspace $ \im (N)$. Hence, given an initial condition $x(0) = z \in \mathbb{R}^n_{\geq 0}$, the trajectory is confined in the affine linear subspace $z+ \im (N)$.
 Additionally, both the positive and non-negative real orthants,
 $\R_{> 0}^n$ and $\R_{\geq 0}^n$, 
  are forward invariant 
with respect to the mass-action ODE system~\eqref{eq:ode}~\cite{volpert}. 
A set is {\em forward invariant} if trajectories that begin in the set remain in the set for all positive time. 
From our discussion, we see that the following sets are also forward invariant:
\begin{align} \label{eq:scc}
\mathcal{P}_{z} = (z + \im(N) ) \cap \R_{\geq 0}^n,\qquad z\in \R^n_{\geq 0}~,
\end{align}
 and these sets
 are called 
 \emph{stoichiometric compatibility classes}. 
 As they are the intersection of an affine linear subspace with the non-negative orthant, they are polyhedra, 
 and
  Sturmfels therefore calls these classes `invariant polyhedra'.
 Equations defining $\mathcal{P}_{z}$ or, more precisely, the affine space $z + \im(N)$ are \emph{conservation laws}. 
 
 \begin{example}[Example~\ref{ex:mc}, continued] \label{ex:McKeithan-polyhedra}
 We saw earlier that the McKeithan network has two linearly independent conservation laws: 
 \begin{equation}\label{eq:mckeithan_cons}
 x_1+x_3+x_4=c_1,\qquad x_2+x_3+x_4=c_2,
 \end{equation}
 for some $c_1,c_2 \geq 0$ (which arise from some $z$ as in \eqref{eq:scc}).  
 One stoichiometric compatibility class of the McKeithan network, namely, the one containing the point $z= \left( \tfrac{1}{2},\tfrac{3}{2},\tfrac{1}{4},\tfrac{1}{4} \right)$, is as follows:
\begin{align*} 
\invtPoly_{z} ~=~ \big\{  { (x_1, x_2, x_3, x_4)} & \in  \mathbb{R}_{\geq0}^{4} \mid  x_1+x_3+x_4=1, ~x_2+x_3+x_4=2 ~\big\}~.
\end{align*}
This set $\invtPoly_{z}$ 
is the triangle depicted here:

\medskip
\begin{center}
\begin{tikzpicture}[scale=0.7]
\node[left] (a) at (-0.1,1) {$(0,1,1,0)$ };
\node[left] (a) at (-0.1,-1) {$(0,1,0,1)$ };
 \node[right] (a) at (3.1,0) {$(1,2,0,0)$};
\draw[fill=black!20!white,line width=0.8pt] (0,1) -- (0,-1) -- (3,0) -- (0,1);
\draw[fill=black] (0,1) circle (2pt);
\draw[fill=black] (0,-1) circle (2pt);
\draw[fill=black] (3,0) circle (2pt);
\end{tikzpicture}
\end{center}
\end{example} 

Recall that 
the boundary of a polyhedron is comprised of {\em faces}.  
For instance, $\invtPoly_{z}$ above has 3 {\em vertices} (0-dimensional faces) and 3  {\em edges} (1-dimensional faces).  Similarly,  
a cube has 8  vertices, 12 edges, and 6 {\em facets} (maximal-dimension proper faces).
Each face $F$ of a stoichiometric compatibility class $ \invtPoly$ is defined by the vanishing of a (possibly non-unique) subset of the variables $x_i$, with index set $W \subseteq \{1,\dots,n\}$, as follows:

\begin{equation}\label{eq:faces}
F= F_W= \invtPoly \cap Z_W,\qquad \textrm{where}\qquad  Z_W := \left\{ x \in \mathbb{R}^{n} ~|~ x_{i}=0 \text{ for all } i \in  {W}  \right\}~.
\end{equation}

For instance, for the invariant polyhedron shown earlier, we label two of the faces according to~\eqref{eq:faces}:
  
\medskip
\begin{center}
\begin{tikzpicture}[scale=0.7]
\node[left] (a) at (-0.1,1) {$(0,1,1,0)$ };
\node[left] (a) at (-0.1,-1) {$(0,1,0,1)$ };
 \node[right] (a) at (3.1,0) {$F_{\{3,4\}}=(1,2,0,0)$};
\draw[fill=black!20!white,line width=0.8pt] (0,1) -- (0,-1) -- (3,0) -- (0,1);
\draw[fill=black] (0,1) circle (2pt);
\draw[fill=black] (0,-1) circle (2pt);
\draw[fill=black] (3,0) circle (2pt);
\node (a) at (1.7,0.9) {$F_{\{4\}}$};
\end{tikzpicture}
\end{center}

 \bigskip
 The final objects of interest are \emph{steady states} (or \emph{equilibrium points}). 
 These are the solutions 
$  x\in \R^n_{\geq 0}$
 to 
 the system $ N v(x) = 0$.  In the case of mass-action kinetics, steady states are the non-negative solutions to a system of polynomials.
 A steady state $x$ is a \emph{positive steady state} if $x \in \R^n_{> 0}$.

\subsection{Questions of interest and overview} 
The mathematical field pertaining to the study of reaction networks -- which was initiated in the 1970s by Feinberg, Horn, Jackson, and others -- focuses broadly on understanding the mathematical properties of the ODE system~\eqref{eq:ode} for \emph{unknown} parameter values. 
The rationale behind this arises from the fact that in many applications, parameter values are difficult 
to estimate, and, moreover, a certain degree of variability among individuals and systems exist. 

The theory  has two broad aims:
    \begin{enumerate}
        \item identify characteristics of the reaction network that impose a certain behavior on the ODE system, mainly for mass-action kinetics, and
        \item develop tools to verify whether and where a property arises for some choice of parameter values.
    \end{enumerate}

\noindent 
Such properties of interest include the number of positive steady states (within a stoichiometric compatibility class), the existence of periodic solutions, the convergence of trajectories to a steady state, the existence of steady states at the boundary of the non-negative orthant, and persistence.  
An ODE system is said to be persistent if no trajectory approaches the boundary. Moreover, one goal in this area is to obtain results on such properties that 
do not rely on simulations or numerical approaches.

In this research area, many mathematical techniques have found application, 
including ideas from applied algebra and polyhedral geometry.
This chapter focuses on selected topics in this direction, specifically those where Sturmfels has contributed. These topics include the algebraic structure and convergence to steady state of complex-balanced systems \cite{TDS}; siphons, persistence, and boundary steady states \cite{ShiuSturmfels}; steady-state invariants via the study of matroids \cite{case-study}; polyhedral methods and mixed volume \cite{case-study}; and convex hulls of trajectories \cite{CKLS}. 
These topics are addressed in what follows, and we highlight Sturmfels's contribution in terms of both specific results and ideas he introduced. 	
	
	\section{Complex-balanced steady states and toric dynamical systems} \label{sec:complex-bal}

	\subsection{The human story}
	In a few places in this chapter, we would like to share some personal stories involving Sturmfels.  We begin with some recollections by Shiu:

{\em As mentioned in the introduction, much progress on Bernd's first article on reaction networks happened at a March 2007 workshop at IMA.   
Before the workshop, Bernd met with me and fellow graduate student Jason Morton.  
Bernd went through the workshop's participant list and told Jason and me which researchers each of us had to meet.  At the top of that list were `Feinberg and Fienberg' -- Jason was to talk with Stephen Fienberg  and 
I was to make sure to meet Martin Feinberg.  
  Funny enough, Bernd is seated between Fienberg and Feinberg in the group photo at the workshop  (see Figure~\ref{fig:IMA}).  

My first night at the workshop, 
Bernd had me meet Alicia Dickenstein for dinner (I was nervous to meet such a prominent mathematician), and then, as mentioned earlier, the three of us, plus Gheorghe Craciun, met early each day.  Gheorghe had been Feinberg's Ph.D.\ student, and was already a prominent researcher in chemical reaction network theory.  One morning, Bernd asked what was known about stability of complex-balanced steady states, and Gheorghe explained that these steady states are known to be local attractors but conjectured to be global attractors as well.  `The global attractor conjecture,' Bernd said, and from then on, the name stuck.  More on this conjecture is in Section~\ref{sec:GAC}.

Another term Bernd coined in the article is (in the title!) `toric dynamical systems' (also known as complex-balanced systems).  In June 2007, Bernd wrote the three of us an email, which began, ``Dear all, this morning I picked up Jeremy Gunawardena's article `Chemical reaction network theory for in-silico biologists', and I simply couldn't stop reading. This author is amazing !!''~\cite{Guna}\footnote{Gunawardena appears again later in another human story (Section~\ref{sec:human-invts}).}.  He continued, ``The main point of his paper is exactly what I was trying express all along: Deficiency zero is just a very special case of [toric dynamical systems], and it is the toric structure that really matters [...]''.  

Then, following a detailed math description, he concluded, ``So, now that I finally understand the literature, I would strongly urge that we use the term `toric dynamical system' (TDS) [...]. This ensures in particular what I was trying to sell to you (Anne) on Wednesday in the car: deficiency zero implies TDS.''  Of course, these ideas had been known for some time (see Section~\ref{sec:cpx-bal}), but the algebraic packaging of the results, including a name that was inviting to algebraic geometers, was important for bringing newcomers into this area.  

As for the car ride mentioned in the email, research meetings with Bernd held on modes of transportation were not uncommon!  Once, he drove me and another Ph.D.\ student (Peter Huggins) from a computational biology retreat near Lake Tahoe back to Berkeley, first with Peter in the front seat for a research chat, and then stopping halfway so Peter and I could switch and then it was my turn to have a math conversation.  Efficient!  }

	\subsection{Steady-state structure } \label{sec:cpx-bal}
For a mass-action system~\eqref{eq:odeYCG}, 
a \emph{complex-balanced steady state} is  a point $x^*\in \R^n_{>0}$ such that
\begin{equation}\label{eq:CB}
A_\k  (x^*)^Y =0.  
\end{equation}
Clearly, by \eqref{eq:odeYCG}, every complex-balanced steady state is 
indeed 
a steady state. 
Complex-balanced steady states were introduced by Horn and Jackson \cite{HornJackson1972} and Feinberg \cite{Feinberg1972} in the 1970s, to account for thermodynamic constraints in a reaction network. 
It follows from \eqref{eq:CB}  that complex-balanced steady states have the following characterization:
for each complex $y$,  the sum of the rates of all reactions with source $y$ equals the sum of the rates of reactions with endpoint $y$. 

Feinberg, Horn, and Jackson recognized that the existence of complex-balanced steady states requires the reaction network to be \emph{weakly reversible}, that is, each connected component must be strongly connected~\cite{Horn1972}. In this case, the image of $C_{\mathcal{G}}$ and $A_\k$ coincide for all $\k\in \R^r_{>0}$  \cite{FeinbergHorn1977}. 

Being weakly reversible, however, is not sufficient for the existence of complex-balanced steady states.
A well-known and celebrated result due to Feinberg is the \emph{deficiency zero theorem} \cite{Feinberg1987}. The \emph{deficiency} of a reaction network is defined as follows:
\[\delta = \dim(\ker (Y) \cap \im(C_\mathcal{G})) . \]
Hence, $\delta \geq 0$.  Also, it can be shown that 
\[
\delta =  m - \ell -s~,
\]
where $\ell$ is the number of connected components of $\mathcal{G}$, and $s$ the rank of the stoichiometric matrix $N$. 
Clearly, if $\delta=0$, then all steady states satisfy \eqref{eq:CB} and hence are complex-balanced.
This is the case for the McKeithan network \eqref{eq:mckeithan}, as $\ell=1$, $s=2$ and hence
\begin{align} \label{eq:mck-def}
 \delta = 3 - 1 - 2=0~.
 \end{align}
Therefore, every positive steady state $x^*$ of the McKeithan network ODEs~\eqref{eq:odeex} is a solution to the system~\eqref{eq:CB}, which translates to:
\begin{align*}
-\k_1 x^*_1x^*_2+ \k_2 x^*_3 + \k_4 x^*_4 =0,\qquad \k_3 x^*_3 - \k_4x^*_4=0 
\end{align*}
(one equation is removed as it is in the linear span of the above two). 
By adding the second equation to the first, we obtain an equivalent system
\begin{align}\label{eq:mckeithan_red}
-\k_1 x^*_1x^*_2+ (\k_2 + \k_3) x^*_3  =0,\qquad \k_3 x^*_3 - \k_4x^*_4=0,
\end{align}
which now consists of binomial equations.  For fixed $\k$, the intersection of the variety 
of~\eqref{eq:mckeithan_red} 
with $\R^4_{>0}$ is the positive part of a toric variety. Hence, it is irreducible and admits a parametrization. 
Such remarkable properties hold for all networks with $\delta=0$, and 
are part of 
why complex-balanced steady states have received 
special mathematical attention.

\medskip
For networks with $\delta>0$, it is no longer true that there are complex-balanced steady states for all parameters $\k$. 
The set $\Sigma$ of parameters that do generate complex-balanced steady states was shown 
in the 1970s 
to be cut out by $\delta$ equations~\cite{Horn1972}. 
Building on that work, Craciun, Dickenstein, Shiu, and Sturmfels further elucidated the algebraic structure of the set $\Sigma$ and also the set of complex-balanced steady states~\cite{TDS} (Theorems~\ref{thm:c-bal-via-T}, \ref{thm:cb}, and~\ref{thm:cb-2} below). 
These results are phrased in the language of toric geometry, as the equations \eqref{eq:CB} (it turns out) define a binomial ideal. This is why, in~\cite{TDS}, an ODE system \eqref{eq:ode} for which all steady states are complex-balanced is called a \emph{toric dynamical system}.  

\smallskip
{\bf Setup.} The setup for the remainder of this subsection is as follows.  
Let $\mathcal G$ be a weakly reversible network with $m$ complexes. 
For each $i\in\{1,\dots,m\}$, 
let $I_i$ denote the subset of $ \{1,\dots,m\}$ corresponding to the complexes in the same connected component as the complex $y_i$. Consider the submatrix $A_{\k,I_i}$ of $A_\k$ with rows and columns in $I_i$, and define $K_i\in \Q[\k]$ as  $(-1)^{\#I_i -1}$ times the minor of $A_{\k,I_i}$ obtained by removing the row and column corresponding to $y_i$. 
Every coefficient of the polynomial $K_i$ is positive.  
Also, 
by \cite[Lemma 5]{TDS}, these polynomials $K_i$ are algebraically independent. 

Next, consider the following ideal obtained by saturation\footnote{Saturation of binomial ideals, especially in the context of algorithms, is an old idea that Bernd has been harnessing for many years; see the work of Bernd with Ho\c{s}ten~\cite{grin}.}:
\[T_{\mathcal{G}} = \langle K_i x^{y_j} - K_j x^{y_i} \mid y_i\rightarrow y_j \in \mathcal{R}   \rangle \colon (x_1\cdots x_n)^\infty 
~ \subseteq~ \mathbb{Q}[x,K] ~\subseteq ~ \mathbb{Q}[x,\k].\]
This ideal gives a characterization of complex-balanced steady states, in the following sense \cite{TDS}.
\begin{theorem}[Complex-balanced steady states] \label{thm:c-bal-via-T}
For $\k^* \in \mathbb{R}^{r}_{>0}$ and $x^* \in \mathbb{R}^n_{\geq 0}$, we have that 
 $x^*$ is  a complex-balanced steady state 
 of the mass-action system of $\mathcal G$ with rate constants $\k^*$ if and only if 
$(x^*,\k^*) $ lies is in the positive part of the toric variety defined by $T_{\mathcal{G}}$, denoted by 
$V_{>0}(T_{\mathcal{G}})$.
\end{theorem}

To determine for which parameter values $\k^*$ there are points of the form $(x^*,\k^*)$ in $V_{>0}(T_{\mathcal{G}})$ (equivalently, which rate constants give rise to complex-balanced steady states), we consider the following elimination ideal:
\[M_\mathcal{G}= T_{\mathcal{G}}  \cap \mathbb{Q}[\k]. 
\] 
The following is one of the main results in~\cite{TDS}.

\begin{theorem}[Rate constants for complex-balancing, part 1]
\label{thm:cb}
The mass-action system of a reaction network $\mathcal G$ with rate constants $\kappa^*\in \R^{r}_{>0}$ admits a complex-balanced steady state if and only if $\k^* \in V_{>0}(M_\mathcal{G})$. Furthermore, the codimension of $M_{\mathcal{G}}$ equals the deficiency $\delta$ of $\mathcal G$. 
\end{theorem}

Next, we illustrate Theorem~\ref{thm:cb} through two examples, one with $\delta=0$ and one with $\delta>0$.

\begin{example}[McKeithan network] \label{ex:mckeithan-def}
Let us find $V_{>0}(M_\mathcal{G})$ for the McKeithan network. We must first compute the ideal $T_\mG$. To this end, as there is one connected component and $m=3$, we obtain the following using the matrix $A_{\kappa}$ shown earlier in Example~\ref{ex:mc}: 
\begin{equation}\label{eq:K} K_1 = (\k_2+\k_3)\k_4,\qquad K_2=\k_1\k_4,\qquad K_3 = \k_1\k_3.
\end{equation}
The following {\tt Macaulay2} code computes the ideal $T_\mG$:
\begin{verbatim}
	R = QQ[x1,x2,x3,x4, K1,K2,K3, MonomialOrder => Eliminate 4];
	I = ideal(K1*x3-K2*x1*x2, K2*x4-K3*x3, K1*x4-K3*x1*x2);
	TG = saturate(I, x1*x2*x3*x4);
	groebnerBasis TG
\end{verbatim}
The output reveals that 
\begin{align*}
T_\mG 
&=  \langle K_1 x_3 - K_2 x_1x_2, ~K_2x_4 - K_3x_3,~K_1x_4 - K_3 x_1x_2 \rangle.
\end{align*}
No polynomials purely in the $K_i$'s appear above as generators of $T_\mG $, and so, as a monomial ordering that eliminates the $x_i$'s was used, 
the elimination ideal $M_\mG$ is the zero ideal. This is consistent with Theorem~\ref{thm:cb}, as $\delta=0$. 
\end{example}

\begin{example}[Extended McKeithan network] \label{ex:extended-mckeithan}
We add the pair of reactions $2X_3\ce{<=>[\k_5][\k_6]} 2X_4$ to the McKeithan network. The resulting {\em Extended McKeithan network} has two connected components ($\ell=2$) and $\delta=1$. 
We let
$y_4=2X_3$ and $y_5=2X_4$, and so we add
 $K_4= \k_6$ and $K_5=\k_5$ to \eqref{eq:K}.
A computation like the one in Example~\ref{ex:mckeithan-def} reveals the following ideal (where the generators form a Gr\"obner basis with respect to a monomial ordering for eliminating $x_1,x_2,x_3,x_4$):
\begin{multline*} T_\mG  =\langle K_{2}^{2} K_{5}-K_{3}^{2} K_{4},~ K_{2} x_{4}-K_{3} x_{3},~ 
K_{4} x_{4}^{2}-K_{5} x_{3}^{2},~ K_{2} x_{1} x_{2}-K_{1} x_{3},~  \\   \qquad 
K_{3} x_{1} x_{2}-K_{1} x_{4},~ K_{3} K_{4} x_{4} -K_{2} K_{5} x_{3} \rangle.
\end{multline*} 
We conclude:
\begin{align} \label{eq:M-ideal}
		M_\mG = \langle K_{2}^{2} K_{5}-K_{3}^{2} K_{4} \rangle ~=~ \langle \k_1^2 \k_4^2 \k_5 - \k_1^2 \k_3^2 \k_6 \rangle. 
\end{align}
Hence, if $\k$ satisfies $ \k_1^2 \k_4^2 \k_5 = \k_1^2 \k_3^2 \k_6 $, 
then there are complex-balanced steady states. Otherwise, there are none.
\end{example}

In Examples~\ref{ex:mckeithan-def} and~\ref{ex:extended-mckeithan}, we determined which rate constants yield complex-balanced steady states by computing the ideals $T_\mG$ and $M_\mG$.  Next, we explain how to avoid having to compute these ideals.  
Namely, Theorem~\ref{thm:cb-2} below gives explicit equations defining $V_{>0}(M_\mathcal{G})$, and they are always binomial in the $K_i$'s.

Reorder the complexes so that complexes in the same connected component are together.  
Let $Y_1,\dots,Y_\ell$ be the submatrices of $Y$ corresponding to connected components. 
Then the \emph{Cayley matrix} is defined as
\begin{align} \label{eq:cayley}
\text{Cay}_\mG(Y)=  \begin{bmatrix}
Y_1 & Y_2 & \dots & Y_\ell \\ 
\mathbf{1} & 0 & \dots & 0 \\
0 & \mathbf{1}  & \dots & 0 \\
\vdots & \vdots & \ddots & \vdots \\
0 & 0 & \dots & \mathbf{1} 
\end{bmatrix} \in \R^{(n+\ell)\times m },
\end{align}
where 
$\mathbf{1}$ is a row vector of $1$'s of the appropriate size. 
With this in place, we can state the next result~\cite{TDS}.

\begin{theorem}[Rate constants for complex-balancing, part 2]
\label{thm:cb-2}
For a reaction network $\mathcal G$, we have that  $\k^* \in V_{>0}(M_\mathcal{G})$  if and only if $K^u=1$ for all $u\in {\rm Ker} ( {\rm Cay}_\mG(Y) )\cap \Z^n$.
 \end{theorem}

\begin{example}[Example~\ref{ex:extended-mckeithan}, continued] \label{ex:mckeithan-extended-part-2}
In the extended McKeithan network, we have 
\[ {\rm Cay}_\mG(Y)=  \begin{bmatrix}
1 & 0 & 0 & 0 & 0   \\ 1 & 0 & 0 & 0 & 0  \\ 0 & 1 & 0 & 2 & 0 \\ 0 & 0 & 1  & 0 & 2 \\ 1 & 1 & 1 & 0 & 0 \\ 0 & 0 & 0 & 1 & 1 \end{bmatrix}, \qquad  
	{\rm Ker} ( {\rm Cay}_\mG(Y) ) = 
	{\rm span} \left\{ \begin{pmatrix} 0\\ 2\\ -2\\ -1\\ 1 \end{pmatrix} \right\}.
	\] 
Hence, $\k\in V_{>0}(M_\mathcal{G})$ if and only if the equation $K_1^0 K_2^2 K_3^{-2} K_4^{-1} K_5=1$ holds, or, equivalently,
\[K_2^2K_5- K_3^2K_4=0.\]
As it should, this binomial equation agrees with the one in~\eqref{eq:M-ideal} we obtained earlier. \end{example}

We end this subsection by expanding on Sturmfels's contribution to two of the ideas mentioned above.  The first concerns the Cayley matrix~\eqref{eq:cayley}, which is so-named because it refers to the Cayley trick in elimination theory \cite{TDS}. Notably, the Cayley trick was first 
developed by Sturmfels \cite{sturmfels_resultant}, and later generalized~\cite{huber_cayley}. 

Another notable contribution by Sturmfels concerns the polynomials $K_i$, which we saw arise from minors of $A_\k$.  
These polynomials admit a graphical interpretation as labels of spanning forests of the digraph $\mG$ labeled with the reaction rate constants. 
This is the 
\emph{Matrix-Tree theorem} of Tutte, which concerns minors of a graph's Laplacian matrix~ \cite{stanley:enumerative,Tutte-matrixtree}. Applications of this theorem appear frequently   in the field of reaction networks, for example in \cite{ConradiShiu, Dickenstein:2011p1112, feliu-wiuf-ptm,feliu-wiuf-crn,saez-gph,fw2013,saez:linear, translated}, and of course in \cite{TDS}.
In fact, the classical King-Altman method in enzyme kinetics follows from the Matrix-Tree theorem \cite{king-altman}. 
 Gunawardena acknowledges having learned about the Matrix-Tree theorem from Sturmfels in \cite{TG}, and Feliu learned about it from Gunawardena's article. So it is probably fair  to say that it was Sturmfels who brought it into the field!

\subsection{Global attractor conjecture} \label{sec:GAC}
A mass-action system that has a complex-balanced steady state is usually called a complex-balanced system, or toric dynamical system. Horn and Jackson proved that for these systems,	 each stoichiometric compatibility class $\invtPoly$, as in~\eqref{eq:scc}, has a unique positive steady state $x^*$, and that $x^*$ is in fact complex-balanced~\cite{HornJackson1972}.  
 Sturmfels gave the name {\em Birch point} to $x^*$, to emphasize the connection to Birch's theorem in (algebraic) statistics\footnote{The existence and uniqueness of positive steady states in toric dynamical systems is essentially equivalent to the existence and uniqueness of maximum likelihood estimators of toric, i.e., log-linear, statistical models. Extensions of Birch's theorem from the study of chemical reaction networks are found in~\cite{GMS2,generalize-birch}.}.	 	 

Next, given a Birch point $x^*$ of a complex-balanced system, the following is a strict Lyapunov function (with respect to the stoichiometric compatibility class $\invtPoly = \invtPoly_{x^*}$): 
\[
 \sum_{i=1}^n \bigl(\, x_i  {\rm log}(x_i) \,
-\, x_i  {\rm log}(x_i^*)  \,-\, x_i \, + x_i^*\bigr)~.
\]
This Lyapunov function ensures that trajectories in $\invtPoly$ that start near the Birch point $x^*$ converge to $x^*$~\cite{HornJackson1972}; in other words, $x^*$ is a local attractor (relative to $\invtPoly$).  The following conjecture states that $x^*$ is also a global attractor. 
	\begin{conjecture}[Global attractor conjecture] \label{conj:GAC} 
	Consider a complex-balanced system (a toric dynamical system).  
	Let $\invtPoly$ be a stoichiometric compatibility class of the system, and let $x^*$ denote the Birch point (the unique positive steady state in $\invtPoly$).  Then every trajectory beginning in the relative interior of $\invtPoly$ converges to $x^*$.
	\end{conjecture}	

As mentioned earlier, the name of this now-famous conjecture was coined by Sturmfels, and this attention revived interest in the conjecture, even though it had been stated much earlier, by Horn, in 1974~\cite{Horn74}.  Another contribution of Sturmfels and co-authors was to reframe what needs to be proven via the language of polytopes, as we explain next.  

The Lyapunov function ensures that every 
trajectory in $\invtPoly$ either converges to $x^*$ or converges to a point on the boundary of $\invtPoly$~\cite{SiegelChen94}.  Such a point on the boundary is necessarily a (boundary) steady state; this fact was shown in a special case by Sontag~\cite{Sontag01}, and the general result follows in a similar manner.  In summary, the global attractor conjecture is equivalent to the following assertion: {\em for a toric dynamical system, no trajectory beginning in the positive orthant converges to a boundary steady state.}
It follows that if $\invtPoly$ has no boundary steady states, then the conjecture holds for
this $\invtPoly$~\cite{SiegelMacLean}.

As noted earlier, Sturmfels often calls stoichiometric compatibility classes `invariant polyhedra', and in~\cite{TDS} (and, independently, by Anderson~\cite{Anderson08}), it is shown that vertices of such polyhedra cannot be points of convergence.  This result set the stage for future progress on the conjecture -- to rule out faces of polyhedra according to their dimension.  Indeed, subsequently, relative-interior points of two types of faces were ruled out: facets (i.e., top-dimension, proper faces)~\cite{AndersonShiu10} and 
``weakly dynamically non-emptyable'' faces~\cite{JohnstonSiegel}.
A related family of ideas 
contributed to a resolution of the conjecture for more networks: those having only one connected component~\cite{Anderson11} and those for which
the stoichiometric compatibility classes 
$\invtPoly$ have dimension at most 3~\cite{CNP,Pantea}.

\begin{figure}[htbp]
\begin{center}
	\begin{align*} 
	\begin{xy}<9mm,0cm>:
	(5,0) 	="y1" *{\bullet} *+!R{0} ; 
	(7,0) 	="y2" *{\bullet} *+!L{2A}  *{}; 
	(5,1)	="y2'" *+!R{B}  *{}; 
	(5,2) 	="y3" *+!R{2B}  *{\bullet}; 
	(6,1) 	="y3'" *{} *+!L{A+B}  *{}; 
	{\ar "y1";"y3'"*{}  }; 		
	{\ar "y2";"y2'"*{}  }; 		
	{\ar "y3";"y3'"*{}  }; 		
	"y1";"y2" **\dir{.};	
	"y2";"y3" **\dir{.};	
	"y3";"y1" **\dir{.};	
	(0,0) 	="A" *{\bullet} *+!R{C} ; 
	(2,1) 	="2A+B"*{\bullet}*+!L{3C+D}  *{}; 
	(0,1) 	="C+D" *{\bullet} *+!L{C+D} ; 
	(-1,2) 	="2D" *{\bullet}*+!R{2D}  *{}; 
	{\ar "A";"2A+B"*{}  }; 		
	{\ar "2A+B";"A"*{}  }; 			
	{\ar "C+D";"2D"*{}  }; 		
	{\ar "2D";"C+D"*{}  }; 		
	"A";"2D" **\dir{.};	
	"2A+B";"2D" **\dir{.};	
	\end{xy}
	\end{align*}
\caption{The network on the left is endotactic, while the one on the right is strongly endotactic.  Both networks are endotactic, because no reaction arrows point out of the Newton polytopes (the triangles) unless they first pass through the respective triangle.  
The network on the right is strongly endotactic, because every proper face of the triangle (that is, each edge and vertex) has a reaction arrow that starts on that face and exits the face.  In contrast, the network on the left is not strongly endotactic: the edge of the triangle from $C$ to $3C+D$ has no reaction arrow that leaves that edge. }
\label{fig:endotactic}
\end{center}
\end{figure}
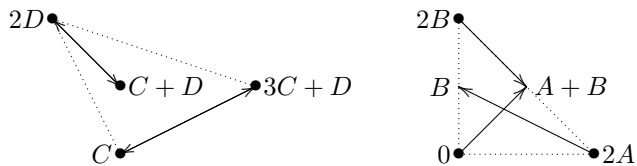

Furthermore, the vocabulary of polytopes led researchers to define `endotactic'~\cite{CNP} and `strongly endotactic' networks~\cite{ProjArg}; see Figure~\ref{fig:endotactic}.  Such networks, roughly speaking, have reaction arrows that point inward, rather than outward (`endo' means `inward', and `tactic' refers to movement as in `chemotaxis').  Geometrically, this means that arrows do not point `out' of the Newton polytope defined by the right-hand sides of the mass-action ODEs; Newton polytopes are discussed in Section~\ref{sec:mss}.  Intuitively, these inward-pointing arrows should guarantee that species concentrations avoid going to 0 or infinity -- that is, trajectories avoid converging to the boundary.  Indeed, for strongly endotactic networks, the global attractor conjecture has been proven~\cite{GMS2}.  This conjecture, however, is still open for endotactic networks.

To summarize, the remaining cases for
Conjecture~\ref{conj:GAC} are toric dynamical systems of dimension at least 4 for which the boundary steady states lie on faces of dimension at least~1 and
codimension at least~2.  See Table~\ref{tab:GAC}.

\begin{table}[ht]
\caption{Properties that guarantee that the global attractor conjecture holds.}
\label{tab:GAC}
\begin{center}
\begin{tabular}{l c}
\hline
Property & Reference \\
\hline
No boundary steady states & \cite{SiegelMacLean} \\ 
No relevant siphons (see Section~\ref{sec:siphon}) & \cite{PetriNetExtended} \\
Strongly endotactic & \cite{GMS2}
\\
Dimension at most $3$ & \cite{Pantea}
\\
Each boundary steady state is a vertex  or & \cite{Anderson08,AndersonShiu10,TDS} \\  \hspace{0.5cm} interior point of a facet 
\\
\hline
\end{tabular}
\end{center}

\end{table}%

The most recent progress on Conjecture~\ref{conj:GAC} is due to Gheorghe Craciun, who proposed an outline of a proof in 2015~\cite{GAC}.  One crucial idea is again to show how the boundary `repels' trajectories, but getting everything to fit together is very subtle.  
As this conjecture was the foremost open problem in the research area, Sturmfels initiated a workshop at San Jos\'e State University in 2016 on Craciun's proposed proof~\cite{gac-siam-news}.  Participants -- including members of Sturmfels's research group -- were assigned parts of the paper to present.  And, at Sturmfels's request, Craciun himself was in attendance but not allowed to speak during the presentations!  
The first installment of Craciun's work is now published~\cite{CraciunSIAGA}, 
and our community eagerly awaits the next one. 

Finally, we note that ideas pertaining to being endotactic have grown outward (`exotactically' perhaps?) into new domains.  
For instance, the concept of `strongly endotactic' was used to analyze stochastic systems 
\cite{a-d-e,a-d-e2, anderson2020tier}. 
Also, results involving strongly endotactic networks (and pertaining to {\em persistence}, a topic of the next section) have been transferred to networks related to some `origin of life' models~\cite{craciun2020autocatalytic}.  Bernd's influence has therefore come full circle: he brought new ideas to the reaction network community, and now ideas developed in this area are in turn proving to be useful in new domains.

	\section{Siphons, persistence, and boundary steady states} \label{sec:siphon}

A~siphon of a reaction network is a subset of species whose absence is
forward-invariant with respect to the dynamics.  
In this section, we explain how siphons are related to persistence and boundary steady states, and elucidate Sturmfels's contribution to exploring the algebraic and computational aspects of siphons.

Let us describe the motivation behind siphons. 
In the McKeithan network~\eqref{eq:mckeithan}, 
all reactions are `off' if and only if either (i) the concentrations of $X_1$, $X_3$, and $X_4$ are all zero, or (ii) the concentrations of $X_2$, $X_3$, and $X_4$ are all zero.  Accordingly, we say that 
$\{ X_1, X_3, X_4 \}$ and $\{ X_2, X_3, X_4 \}$ 
(or, more accurately, $\{1,3,4\}$ and $\{2,3,4\}$) 
are (minimal) {\em siphons}. 

More precisely, a {\em siphon} of a network is a non-empty subset $W \subseteq \{1,\dots,n\}$ of the set of species indices, such that if the product complex of a reaction contains an element of $W$, then so does 
  the reactant complex. For example, if  $X_1+X_2 \rightarrow X_3$ is a reaction of a network, then 
 $3 \in W$
 implies that
 $1 \in W$
 or 
 $2 \in W$.
 A siphon of a network $\mathcal G$ is {\em minimal} if it is minimal with respect to inclusion among all siphons of $\mathcal G$.

The connection between siphons, faces of $\invtPoly$, and boundary steady states is the following result~\cite{Anderson08,PetriNetExtended}, where we
recall that the face $F_W$ induced by $W$ was defined in \eqref{eq:faces}:

\begin{proposition}
If there is a steady state in the relative interior of some face $F_W$ of $\invtPoly$, then $W$ is a siphon.
\end{proposition}

This proposition -- and our interest in the global attractor conjecture and persistence of mass-action systems  -- motivate the need for efficient tools for computing siphons.
A mass-action system is {\em persistent} if 
no trajectory beginning in the relative interior of a stoichoimetric compatibility class $\invtPoly$ has an accumulation point, or {\em $\omega$-limit point}, 
on the boundary of $\invtPoly$. 
In this context, Shiu and Sturmfels proved the following two results~\cite{ShiuSturmfels}:

\begin{proposition}
 If a reaction network is strongly connected, then its minimal siphons correspond exactly to minimal associated primes of the ideal of complex monomials $x^{y_{i}}$.
 \end{proposition}
 
For instance,  the ideal of complex monomials for the McKeithan network~\eqref{eq:mckeithan} is $\langle x_1x_2,x_3,x_4 \rangle$ in $\mathbb{Q}[x_1,x_2,x_3,x_4]$, and its minimal associated primes are $\langle x_1,x_3,x_4 \rangle$ and $\langle x_2,x_3,x_4 \rangle$.  These primes correspond exactly to the minimal siphons found earlier.

To state the more general result, which allows for networks that are not strongly connected, consider the following ideal in the quotient ring $\mathbb{Q}[x_1, \dots, x_n]/ \langle x_1 x_2 \cdots x_n \rangle$:
\[
\idealI_{\mathcal G} ~=~
\bigl\langle \, x^{y_{i}} \cdot (x^{y_{j}} - x^{y_{i}})  ~:~ {y_{i}} \rightarrow {y_{j}} \text{ is a reaction of } 
\mathcal{G} \,\bigr\rangle~.
\]

\begin{theorem}[Minimal siphons]\label{thm:siphon}  
Let $\mathcal G$ be a reaction network.  For each minimal prime~$\idealP$ of $\idealI_{\mathcal G}$, consider the set:
\begin{align} \label{eq:siphon-set}
    \{i \in \{1,\dots,n\} \mid x_i \in \idealP \}.
\end{align}
The inclusion-minimal sets of the form~\eqref{eq:siphon-set} are precisely the minimal siphons of~$\mathcal G$. 
\end{theorem}

\begin{example}
Consider the following network (its biological significance will be explained in Section~\ref{sec:idea-invariant}):
\begin{equation}\label{eq:G1-after-relabeling}
\begin{aligned}
& X_1 +X_3  \to X_6   \to X_1+X_4  \to X_7   \to X_1 + X_5 
\\
& X_2 +X_5  \to X_8   \to X_2+X_4  \to X_9   \to X_2 + X_3. \end{aligned}
\end{equation}
As in~\cite{ShiuSturmfels}, we use {\tt Macaulay 2}~\cite{M2} and, in particular, the command {\tt decompose} to obtain minimal primes:
\begin{verbatim}
R = QQ[x_1 .. x_9];
I = ideal (x_1*x_3*(x_6-x_1*x_3),  x_6*(x_1*x_4-x_6),  
 x_1*x_4*(x_7-x_1*x_4), x_7*(x_1*x_5-x_7), x_2*x_5*(x_8-x_2*x_5),  
 x_8*(x_2*x_4-x_8), x_2*x_4*(x_9-x_2*x_4),  x_9*(x_2*x_3-x_9));
decompose (I + ideal product gens R)
\end{verbatim}
The output reveals that there are 12 minimal primes:
\begin{align*}
& \langle x_{1},\:x_{2},\:x_{6},\:x_{7},\:x_{8},\:x_{9}\rangle,~
 \langle x_{3},\:x_{4},\:x_{5},\:x_{6},\:x_{7},\:x_{8},\:x_{9} \rangle,~
 \langle x_{1},\:x_{4},\:x_{5},\:x_{6},\:x_{7},\:x_{8},\:x_{9} \rangle, \\
&
  \langle x_{2},\:x_{3},\:x_{4},\:x_{6},\:x_{7},\:x_{8},\:x_{9} \rangle,~ 
  \langle x_{2},\:x_{3},\:x_{4},\:x_{6},\:x_{8},\:x_{9},\: x_1 x_5 - x_7\rangle,~ \\
  & 
  \langle x_{2},\:x_{3},\:x_{6},\:x_{8},\:x_{9},\:x_4 - x_5,\: x_1 x_5 - x_7 \rangle, ~
   \langle x_2,\: x_3, \: x_8,\: x_9,\: x_6 - x_7,\: x_4 - x_5, \: x_1 x_5 - x_7 \rangle,~
\\ &
    \langle 
    x_2,\: x_8,\: x_9,\:
    x_3 - x_5,\:
    x_4 - x_5,\:
    x_6 - x_7, \:
     x_1 x_5 - x_7 \rangle,~
     \langle x_{1},\:x_{4},\:x_{5},\:x_{6},\:x_{7},\:x_{8},\: x_2  x_3 - x_9 \rangle,~
  \\
  &
      \langle x_{1},\:x_{5},\:x_{6},\:x_{7},\:x_{8}, \: x_3 - x_4, \: x_2 x_4 - x_9 \rangle,~ 
      \langle x_{1},\:x_{5},\:x_{6},\:x_{7},\: x_3 - x_4,\: x_8 - x_9,\:  x_2 x_4 - x_9 \rangle,~ 
  \\
  &      \langle x_1,\: x_6,\: x_7,\: 
		 x_3 - x_5,\:   
		x_4 - x_5,\: 
  		x_8 - x_9,\:
		 x_2 x_5 - x_9
		 \rangle~.
\end{align*} 
The last prime ideal contains $x_1,x_6,x_7$.
Additionally, for all other prime ideals, the set of all monomials of the form $x_i$ appearing in the ideal is {\em not} a proper subset of $\{x_1, x_6, x_7\}$. 
Thus, $\{1,6,7\}$ is a minimal prime.  Similarly, there are two more minimal primes: $\{2,8,9\}$ and $\{3,4,5,6,7,8,9\}$.
\end{example}

The next result involves a certain square-free monomial ideal $\mathfrak B$, which is constructed from the stoichiometric subspace (see~\cite{ShiuSturmfels} for details).  
The result concerns relevant siphons; a siphon $W$ is {\em relevant} if it defines a non-empty face of some stoichiometric compatibility class $\invtPoly$ (i.e., the face $F_W$ of $\invtPoly$ is non-empty).  Shiu and Sturmfels proved the following result, which extends a result of Angeli, De Leenheer, and Sontag~\cite{PetriNetExtended}.

\begin{theorem} \label{thm:relevant-siphon} For a reaction network $\mathcal G$, the following are equivalent:
\begin{enumerate}
	\item For every siphon  $W$ of $\mathcal G$, 
	there exists a non-negative, non-zero\footnote{The fact that `non-zero' is missing from the statement of this result in~\cite{ShiuSturmfels} was pointed out by Gilles Gnacadja (Amgen), whose connection to our community was initiated by an email from Sturmfels in 2007 and then a visit to Berkeley.}, linear conservation law  
 whose support is a subset of $W$.		
	\item  
	$\mathcal G$ has no relevant siphons.  
	\item ${\rm Sat}(\idealI_{\mathcal G},\idealB) = \langle 1 \rangle$.
\end{enumerate}
If the above (equivalent) conditions hold, then there are no boundary steady states in any stoichiometric compatibility class, and in particular the global attractor conjecture holds for every toric dynamical system arising from $\mathcal G$.
\end{theorem}

\begin{example} \label{ex:mckeithan-relevant-siphons}
Recall that the minimal siphons of the McKeithan network are 
$\{1,3,4\}$ and $\{2,3,4\}$.
 Each of these contains the support of a non-negative, non-zero, linear conservation law, as in \eqref{eq:mckeithan_cons}, namely, $x_1+x_3+x_4 - c_1 =0$ and $x_2+x_3+x_4- c_2=0$, respectively.  So, by Theorem~\ref{thm:relevant-siphon},  no boundary steady states exist in any stoichiometric compatibility class.
\end{example}

We note that Sturmfels, motivated in part by computations involving siphons, directly encouraged Thomas Kahle to develop 
theoretical and computation tools for binomial primary decomposition~\cite{kahle2010decompositions,kahle-miller} (which of course built on Sturmfels's foundational work with Eisenbud~\cite{binomial-ideals}).  
Much of Kahle's work was in collaboration with Ezra Miller, and the two of them were first introduced by one of Bernd's famous `I think you should talk' emails.  

Indeed, Bernd truly succeeded in bringing algebraically inclined people into the reaction network community:  Kahle later collaborated with Carsten Conradi (who first learned about algebraic aspects of this area from Karin Gatermann) on detecting whether an ideal (especially a steady-state ideal) is binomial~\cite{conradi-kahle}.  Subsequently, Kahle and Conradi,
 with their jointly advised Ph.D.\ student, Alexandru Iosif, analyzed multistationarity for systems having a monomial steady-state parametrization~\cite{conradi2019multistationarity}.  Also, Ezra Miller went on to prove, with Manoj Gopalkrishnan and Anne Shiu, the strongly endotactic case of the global attractor conjecture~\cite{GMS2,ProjArg} (recall Section~\ref{sec:GAC}).

We end this section by describing one more aspect of Sturmfels's work that pertains to boundary steady states. 
When studying boundary steady states, we often want to know the possible configurations of boundary faces $F_W$ in the various stoichiometric compatibility classes.     
These classes~\eqref{eq:scc} are subsets of $\mathbb{R}^n_{\geq 0}$ and can be defined by means of an equation $Z x=c$, where $Z$ is a ${d\times n}$-matrix.  
In many applications, the image of $\R^n_{\geq 0}$ under $Z$ is $\R^d_{\geq 0}$, and so the stoichiometric compatibility classes are in bijection with 
 $ \R^d_{\geq 0}$. 
In this setting, Sturmfels and collaborators introduced the {\em chamber complex} of $Z$.  This complex is the partition of  $\R^d_{\geq 0}$ into equivalence classes, where $c,c'\in \R^d_{\geq 0}$ are equivalent if their associated stoichiometric compatibility classes (which are polyhedra) have the same normal fan~\cite{case-study} and in particular are combinatorially isomorphic, with the same arrangement of boundary faces.  
The authors of \cite{case-study} give a description of the chamber complex for a specific biological model: the Wnt shuttle network discussed later in Section~\ref{sec:MV}. 
Coming full circle, a related approach to studying this chamber complex, which involves the ideal $\mathfrak B$ mentioned above, was previously described by Sturmfels in the context of siphons~\cite{ShiuSturmfels}. 
		
\section{Invariants, matroids, and model selection} \label{sec:matroid}
We now leave complex-balanced steady states and boundary steady states to enter a new topic, namely how steady-state invariants can be used for model discrimination (Section~\ref{sec:idea-invariant}) and then discuss how Sturmfels introduced the idea of using matroids to help with this approach (Sections~\ref{sec:human-invts} and~\ref{sec:alg-matroid}).

\subsection{The idea behind invariants} \label{sec:idea-invariant}
In 2007, 
Jeremy Gunawardena  wrote an influential article on how algebraic invariants can be used for non-parametric model selection \cite{G-distributivity} (see also \cite{MG-multisite}). The idea can be formulated as follows. 
Given a reaction network $\mathcal{G}$, consider the \emph{steady-state ideal}~$I$, which is the ideal in $\C(\k)[x]$ generated by the polynomials on the right-hand side of the mass-action ODEs \eqref{eq:odeYCG}.
If $x^*$ is a steady state of a mass-action system defined by $\mathcal{G}$ and a choice of rate constants $\k^*$, 
then essentially every polynomial in $I$ vanishes when evaluated at  $x^*$ and specialized at~$\k^*$.  Here, `essentially'
means that we avoid polynomials in which some coefficient has a zero denominator when specialized at~$\k^*$. 
In this context, an element of $I$ is called a \emph{(steady-state) invariant}. 

Model selection aims at determining which of several competing models best represents a mechanism under study. 
Gunawardena's crucial insight is that by comparing the monomials appearing in the steady-state invariants of two competing models (networks), 
we can use experimental data at steady state 
to 
reject a model or support its validity,  \emph{without needing to know (or estimate) parameter values}. 

To illustrate the idea, we present a slight modification of the models Gunawardena analyzed in \cite{G-distributivity}. We consider a substrate $S$ that exists in three states $S_0,S_1,S_2$. Such states arise for example when $S$ admits the addition of a phosphate group (a process called \emph{phosphorylation}) at two different sites: $S_0$ has no phosphate groups attached, $S_1$ has one phosphate, and $S_2$ has two. 
Phosphorylation is mediated by an enzyme $E$, and the reverse process, called {\em dephosphorylation}, is mediated by the enzyme $F$. One possible mechanism is the following (called \emph{distributive}), where an independent encounter of $S$ with the enzyme $E$ or $F$ is required for  the phosphorylation or dephosphorylation of each site:
\begin{align}\label{eq:G1}
\mG_1\colon\qquad  & E +S_0  \ce{->[\k_1]} Y_1   \ce{ ->[\k_2]} E+S_1  \ce{->[\k_3]} Y_2   \ce{->[\k_4]} E + S_2 \\
& F +S_2  \ce{->[\k_5]} Y_3   \ce{ ->[\k_6]} F+S_1  \ce{->[\k_7]} Y_4   \ce{->[\k_8]} F + S_0. \nonumber
\end{align}
An alternative mechanism (called {\em processive}) allows an encounter with the enzyme $E$ to result in the phosphorylation of one or two sites\footnote{The processive mechanism shown here differs slightly from some in the literature, e.g.~\cite{ConradiShiu}.}. This gives the following {\em mixed-mechanism} reaction network (notice that the dephosphorylation mechanism remains distributive):
\begin{align*}
\mG_2\colon \qquad & E +S_0  \ce{->[\k_1]} Y_1   \ce{ ->[\k_2]} E+S_1  \qquad  Y_1   \ce{->[\k_4]} E + S_2 \\
& F +S_2  \ce{->[\k_5]} Y_3   \ce{ ->[\k_6]} F+S_1  \ce{->[\k_7]} Y_4   \ce{->[\k_8]} F + S_0.
\end{align*}
We rename the species as $X_1=E, X_2=F, X_3=S_0, X_4=S_1, X_5= S_2, X_6=Y_1,X_7=Y_2,X_8=Y_3,X_9=Y_4$.  
The relabeled network $\mG_1$ is the network~\eqref{eq:G1-after-relabeling} from earlier.
Then, using mass-action kinetics, we form the steady-state ideals $I_1$ and $I_2$ of (respectively) $\mG_1$ and $\mG_2$:
\begin{align*}
I_1 &= \langle \k_{3} x_{1} x_{4}-\k_{4} x_{7},~ 
-\k_{3} x_{1} x_{4}-\k_{7}  x_{2}x_{4}+\k_{2} x_{6}+\k_{6} x_{8},~ 
-\k_{5}  x_{2}x_{5}+\k_{4} x_{7},~ \\ & \qquad \k_{1} x_{1} x_{3}-\k_{2} x_{6},~ 
\k_{5}  x_{2}x_{5}-\k_{6} x_{8},~ \k_{7}  x_{2}x_{4}-\k_{8} x_{9}
 \rangle, \\
I_2&=  \langle -\k_{7} x_{2} x_{4}+\k_{2} x_{6}+\k_{6} x_{8},~ 
-\k_{5}  x_{2}x_{5}+\k_{4} x_{6},~ 
\k_{1} x_{1} x_{3}-\k_{2} x_{6}-\k_{4} x_{6},~  \\ & \qquad
\k_{5}  x_{2}x_{5}-\k_{6} x_{8},~  \k_{7}  x_{2}x_{4}-\k_{8} x_{9} \rangle.
\end{align*}

In the networks $\mG_1$ and $\mG_2$, the species $Y_i$ are called {\em intermediate species}, and their concentrations are typically difficult to measure experimentally.  We therefore eliminate the corresponding variables from the steady-state ideals, as follows: 
\begin{align}
	\label{eq:elim-1}
	I_1 \cap \C(\k)[x_1,\dots,x_5] & = \langle  \k_1\k_5x_1x_3x_5 - \k_3\k_7 x_1 x_4^2\rangle,  \\
		\label{eq:elim-2}
	I_2 \cap \C(\k)[x_1,\dots,x_5] & = \langle -\k_4\k_{7} x_{2} x_{4} + \k_5(\k_2+\k_4)x_2x_5,~ \k_1x_1x_3 - \k_7 x_2 x_4 \rangle. 
\end{align}

In these elimination ideals, we focus on $x_3,x_4,x_5$ (corresponding to species $S_0, S_1, S_2$).  
By factoring $x_1$ from the generator in~\eqref{eq:elim-1}\footnote{We could instead avoid this factoring by first saturating the steady-state ideal with respect to $x_1\cdots x_9$.},
we find that, for all parameter values, at steady state (specifically, at positive steady states), 
\begin{align} \label{eq:network-1-line}
x_3x_5 \text{ and } x_4^2 \text{ lie on a line}. 
\end{align}
However, the same does not hold for $\mG_2$. 
For this second network, we factor $x_2$ from the first polynomial, and obtain that 
at steady state, 
\begin{align} \label{eq:network-2-line}
x_4 \text{ and } x_5 \text{ lie on a line}. 
\end{align}

The implications for model selection are as follows.  
By obtaining several sets of experimental data at steady state
for $x_3,x_4,x_5$ (arising from distinct initial conditions $x(0)$),  
we can assess which of the pairs $(x_4,x_5)$
 or $(x_3x_5,x_4^2)$  is more likely to have come from a line. 
  We can then determine whether 
 $\mG_1$ or $\mG_2$ is a more accurate model, or perhaps both models need to be rejected.
 
 Of course this is an idealistic approach, because experimental data contain errors, which are magnified when data values are multiplied (via monomials like $x_3 x_5$), and the model itself is an approximation of reality. Accordingly, this approach requires a proper statistical treatment. 
 Some ideas in this direction were given in \cite{harrington-model}, 
 and related approaches to 
non-parametric model selection have been employed for example in \cite{autokinase}.

Despite these challenges, Gunawardena's basic idea is very appealing, as it circumvents the need to know or estimate the parameters of a candidate model. 
Indeed, identifying suitable invariants -- or at least the subsets of variables that are related by a polynomial in the steady-state ideal (more on this topic is in the next subsections) -- can help with experimental design.  More precisely, researchers can construct experiments to measure relevant subsets of variables so that polynomial relations can be tested.

  \subsection{The human story} \label{sec:human-invts}
In March 2013, Dickenstein, Gunawardena, and Shiu organized a workshop,  
`Mathematical Problems Arising from Biochemical Reaction Networks', 
at the American Institute for Mathematics (AIM). 
Sturmfels came to visit, and in the afternoon, a big audience gathered around him while he was standing at a whiteboard talking about \emph{matroids} and their relation to the problem of steady-state invariants. The word `matroid' was foreign to most of the workshop participants, but the idea stuck. Some of us went back home and started reading on the topic. 

Heather Harrington  was present that day. Subsequently, 
she and collaborators applied ideas Sturmfels initiated -- using matroids for finding steady-state invariants --
to the problem of model selection for an important biological process, the Wnt signaling pathway~\cite{MacLean-Wnt}, 
and later explored this topic further with Sturmfels and others~\cite{case-study}.  
These works constituted part of the thesis of Sturmfels's PhD student Zvi Rosen, 
which was entitled `Algebraic matroids in applications'~\cite{rosen_thesis}.

  \subsection{Algebraic matroids } \label{sec:alg-matroid}
So, what is an algebraic matroid? Or well... what is a matroid?  As it is not the aim of this chapter to give a
 proper  introduction to matroids, we will simply say that it is a pair of objects $(X,\mathcal{I})$ such that $X$ is a finite set and 
 $\mathcal{I}$ is a set of subsets of $X$.  
 The pair must satisfy some independence axioms, 
 and so the elements of $\mathcal{I}$ are called {\em independent sets}~\cite{oxley_matroids}.
 The basic example of a matroid (a linear matroid) comes from a finite set of vectors $X=\{v_1,\dots,v_n\}$ from a vector space: here $\mathcal{I}$ consists of all  linearly independent subsets of $X$.
 
 Matroids have the following {\em downward-closed} property: Every subset of an independent set is an independent set 
 (if $U \subseteq V \in \mathcal{I}$, then $U \in \mathcal{I}$).
It thus makes sense to define the \emph{circuits} of the matroid as the minimal sets (with respect to inclusion) not contained in $\mathcal{I}$. 

Algebraic matroids arise in the context of irreducible algebraic varieties, and independent sets consist of variables that are \emph{not} the support of any polynomial in the defining prime ideal. 
Specifically, consider a prime ideal $I$ in a polynomial ring $K[x]$ where $x=(x_1,\dots,x_n)$, and let $ \mathcal{I}$ denote the following subset of the power set of $\{x_1,\dots, x_n\}$:
\[
\mathcal{I} = \big\{ C \subseteq \{x_1,\dots,x_n\} \mid I \cap K[C]=\{0\} \big\} . \]
Hence, circuits are the minimal subsets of variables for which there is a polynomial in the ideal involving only the variables in the subset. Additionally, for a circuit $C$, the ideal $I \cap K[C]$ is a principal ideal 
and the generator involves all variables  in $C$. Several strategies can be employed to find the algebraic matroid \cite{rosen_thesis}. 
 
In the context of reaction networks, we are interested in algebraic matroids arising from the defining prime ideal $I$ (in $\Q(\k)[x]$) of 
some irreducible component of the {\em steady-state variety} (that is, the zero set defined by the right-hand sides of the ODEs).  When the set of positive steady states admits a rational parametrization (as in the examples below), we will instead consider the prime ideal $I$ of the closure of the image of the parametrization.

Let us find the circuits of the algebraic matroids of the networks $\mG_1$ and $\mG_2$ above.  For the distributive mechanism $\mG_1$, from~\eqref{eq:G1}, we use symbolic software 
to solve the equations obtained by setting to zero the right-hand sides of the ODEs for $i=4,\dots,9$, obtaining: 
\begin{align*}
x_{4} &= \tfrac{ \k_{1} }{ \k_{7} }\, x_{1} x_2^{-1} x_{3} , & 
 x_{5} &= \tfrac{ \k_{1} \k_{3}}{\k_{5}\k_{7} } \, x_{1}^{2}  x_{2}^{-2}  x_{3}  , &  x_{6} &= 
\tfrac{ \k_{1}}{\k_{2}}\, x_{1} x_{3}, \\
  x_{7} &= 
\tfrac{ \k_{1} \k_{3}}{ \k_{7} \k_{4}}\,  x_{1}^{2} x_{2}^{-1}  x_{3}, &  x_{8} &= 
\tfrac{\k_{1} \k_{3}}{ \k_{6} \k_{7} }\, x_{1}^{2} x_{2}^{-1} x_{3}  , &  x_{9} &= 
\tfrac{ \k_{1}}{\k_{8}}\, x_{1} x_{3}.
\end{align*}
This yields a {\em steady-state parametrization}: 
for a given $\kappa$, 
the set of positive steady states is the image of the map $\mathbb{R}^3_{>0} \to \mathbb{R}^9_{>0}$ sending $(x_1,x_2,x_3)$ to $(x_1,\dots, x_9)$, where $x_4,\dots, x_9$ come from the above expressions.  
Using implicitization \cite[\S 3.3]{cox-little-oshea}, we find that the corresponding prime ideal is: 
\begin{multline*} 
\widetilde{I}_1 = \big\langle \k_{4} x_{7}-\k_{6} x_{8}, \k_{2} x_{6}-\k_{8} x_{9}, 
-\k_{5} \k_{8} x_{5} x_{9}+\k_{7}  \k_{6} x_{4} x_{8}, 
\k_{1} \k_{6} x_{3} x_{8}-\k_{3} \k_{8} x_{4} x_{9}, \\ 
\k_{1} \k_{5} x_{3}x_{5}  -\k_{3} \k_{7} x_{4}^{2},  
 \k_{5} x_{2}x_{5}-\k_{6} x_{8}, \k_{7} x_{2} x_{4}-\k_{8} x_{9}, 
\k_{3} \k_{8} x_{1} x_{9}-\k_{6}\k_{7}  x_{2} x_{8},  \\
 \k_{3}x_{1} x_{4}-\k_{6} x_{8},  \k_{1} x_{1}x_{3}-\k_{8} x_{9} \big\rangle.
\end{multline*}
Using Gr\"obner bases and elimination, we find that the minimal sets of variables $\{x_{i_1},\dots,x_{i_q}\}$ for which there is a polynomial in $\widetilde{I}_1\cap \Q(\k)[x_{i_1},\dots,x_{i_q}]$ in those variables are the following:
{\small \begin{align*}  q=2\colon\quad  & \{6, 9\}, \{7, 8\}, 
\\ 
q=3 \colon \quad & \{1, 3, 6\}, \{1, 3, 9\}, \{1, 4, 7\}, \{1, 4, 8\}, \{2, 4, 6\}, \{2, 4, 9\}, \{2, 5, 7\}, \{2, 5, 8\}, \{3, 4, 5\},  \\ 
q=4 \colon\quad & \{1, 2, 3, 4\}, \{1, 2, 3, 5\}, \{1, 2, 3, 7\}, \{1, 2, 3, 8\}, \{1, 2, 4, 5\}, \{1, 2, 5, 6\}, \{1, 2, 5, 9\},   \\ &  \{1, 2, 6, 7\}, 
 \{1, 2, 6, 8\}, \{1, 2, 7, 9\}, \{1, 2, 8, 9\}, \{1, 3, 5, 7\}, \{1, 3, 5, 8\}, \{1, 4, 5, 6\},\\ & 
 \{1, 4, 5, 9\}, \{1, 5, 6, 7\},  \{1, 5, 6, 8\}, \{1, 5, 7, 9\}, \{1, 5, 8, 9\}, \{2, 3, 4, 7\}, \{2, 3, 4, 8\},\\ &
  \{2, 3, 5, 6\}, \{2, 3, 5, 9\}, \{2, 3, 6, 7\},  \{2, 3, 6, 8\}, \{2, 3, 7, 9\}, \{2
, 3, 8, 9\}, \{3, 4, 6, 7\},  \\ & \{3, 4, 6, 8\}, \{3, 4, 7, 9\}, \{3, 4, 8
, 9\}, \{3, 5, 6, 7\}, \{3, 5, 6, 8\}, \{3, 5, 7, 9\}, \{3, 5, 8, 9\}
, \\ &  \{4, 5, 6, 7\}, \{4, 5, 6, 8\}, \{4, 5, 7, 9\}, \{4, 5, 8, 9\}.
\end{align*}}%
These sets are the circuits of the matroid. We recognize the circuit $\{3, 4, 5\}$ 
as coming from a steady-state invariant we found earlier (recall~\eqref{eq:network-1-line}).
As we started with 
an irreducible component intersecting $\R^9_{>0}$, factorization of $x_1$ is not necessary in this approach.

We repeat with the mixed-mechanism network $\mG_2$ and find a steady-state parametrization given by:
\begin{align*}
x_{4} &= \tfrac{ \k_{1}}{ \k_{7}}\,  x_{1} x_{2}^{-1} x_{3}, &  x_{5} &= 
\tfrac{\k_{1} \k_{4}}{  (\k_{2}+\k_{4})\k_{5}} \,  x_{1}x_{2}^{-1}x_3 , &  
x_{6} &= \tfrac{ \k_{1}}{\k_{2}+\k_{4}}\, x_{1} x_{3},   \\   x_{8} &= 
\tfrac{\k_{1} \k_{4}}{(\k_{2}+\k_{4}) \k_{6}}\,  x_{1} x_{3}, &  x_{9}
 &= \tfrac{ \k_{1}}{\k_{8}}\, x_{1} x_{3},
\end{align*}
where we recall that $x_7$ is missing simply because $\mG_2$ does not contain the species $X_7$.  
The corresponding prime ideal is as follows:
\begin{multline*}
\widetilde{I}_2=\big\langle (\k_{2} \k_{6}+\k_{6} \k_{4}) x_{8}-\k_{4} \k_{8} x_{9}, 
x_{6} (\k_{2}+\k_{4})-\k_{8} x_{9}, 
(-\k_{2} \k_{5}-\k_{4} \k_{5}) x_{5}+\k_{4} \k_{7} x_{4},  \\
x_{2} x_{5} (\k_{2} \k_{5}+\k_{4} \k_{5})-\k_{4} \k_{8} x_{9}, 
x_{3} x_{1} \k_{1}-\k_{8} x_{9}^{2}\big\rangle. 
\end{multline*}
The circuits of the corresponding algebraic matroid are the following: 
{\small \begin{align*}  q=2\colon\quad  & \{4, 5\}, \{6, 8\}, \{6, 9\}, \{8, 9\}, \\
q=3\colon\quad  &  \{1, 3, 6\}, \{1, 3, 8\}, \{1, 3, 9\}, \{2, 4, 6\}, \{2, 4, 8\}, \{2, 4, 9\}, \{2, 5, 6\}, \{2, 5, 
8\}, \{2, 5, 9\}, \\ 
q=4\colon\quad  &  \{1, 2, 3, 4\}, \{1, 2, 3, 5\}.
\end{align*}}
As expected, $\{4, 5\}$ is a circuit (recall the steady-state invariant coming from~\eqref{eq:network-2-line}).

By comparing the two algebraic matroids, we see that the set $\{1, 2, 3, 4\}$ is a circuit of both. 
So, we consider the elimination ideals $\widetilde{I}_j\cap \Q(\k)[x_1,x_2,x_3,x_4]$, for $j=1,2$, and 
examine the respective generators:
\begin{align*}
\k_1x_1x_3 - \k_7x_2x_4 \text{ (for $\mG_1$)} \quad \text{and} \quad \k_1\k_8x_1x_3-\k_7^2x_2^2x_4^2 \text{ (for $\mG_2$)}. 
\end{align*}  
The monomials in these two polynomials differ, so we conclude that steady-state measurements of $X_1,X_2,X_3,X_4$ could lead to model discrimination.

\subsection{Steady-state parametrizations} \label{sec:steady-state-params}
We end this section by discussing steady-state parametrizations and their relation to matroids.  

In the previous subsection, we saw steady-state parametrizations for $\mG_1$ and $\mG_2$.  
These parametrizations consist of monomials in the $x_i$'s, so the positive part of the steady-state variety 
is toric (recall that the ideals $\widetilde{I}_1$ and $\widetilde{I}_2$ are binomial).  
As another example, we revisit the McKeithan network. The steady-state equations~\eqref{eq:mckeithan_red} are linear in $x_3$ and $x_4$, so we solve for those variables to obtain a steady-state parametrization consisting of monomials:
\begin{equation}\label{eq:mckeithanpar}
x_3 = \tfrac{\k_1}{\k_2+\k_3} x_1x_2, \qquad x_4=\tfrac{\k_1\k_3}{(\k_2+\k_3)\k_4} x_1x_2.
\end{equation}

Of course there is no guarantee that a steady-state parametrization exists, but it happens surprisingly often in networks arising in biological applications.  This is because polynomial systems arising from reaction networks typically have  only linear or quadratic terms, as they represent interactions of proteins. 
This is often enough to enable the \emph{linear elimination} of variables at steady state.  Furthemore, under favorable conditions, the Matrix-Tree theorem (our old friend!) guarantees that the expressions found after elimination are rational functions with positive coefficients -- and hence define the positive part of the steady-state variety; see e.g.~\cite{feliu-wiuf-ptm,feliu-wiuf-crn,saez:linear,linear-framework,TG}.

Steady-state parametrizations can greatly aid in the analysis of a reaction system~\cite{cs-survey}, so we are interested in more strategies for finding them.  One such strategy uses tools from toric geometry, in particular when  the set of positive steady states define a toric variety, e.g. \cite{TSS, J-M-P,messi}. 
Another strategy, employed by Sturmfels and collaborators, uses the 
bases of degree $1$ of the algebraic matroid~\cite[\S 5]{case-study}.

We end this section by noting a final application of matroids to reaction networks: 
Sturmfels and co-authors used matroids to assess {\em identifiability} of parameters, that is, the question of whether parameters -- here, rate constants -- can be recovered from data~\cite[\S 7]{case-study}. 
 
\section{Polyhedral methods and multistationarity} \label{sec:mss}
As the reader might have noticed so far, algebraic and polyhedral geometric tools go hand in hand in the study of steady states of reaction networks. 
We will see more instances of this phenomenon in this section (for assessing whether there are multiple steady states) and in Section~\ref{sec:counting} (for counting steady states). 

\subsection{Multistationarity and injectivity. } \label{sec:mss-injective}
Here we consider the important problem of deciding, for a given network, 
 whether there exist parameter values for which there is more than one positive steady state in some stoichiometric compatibility class. 
When that happens, we say that the network is {\em multistationary}.  Multistationarity has implications in cellular decision-making, and therefore answering this question has been the focus of much research. 

Concretely, we consider a network $\mathcal G$, and let 
$Z\in \R^{d\times n}$ be a matrix such that the stoichiometric compatibility classes are of the form $Zx=c$ for $c\in \R^d$ (so $d=n-s$, and $s$ is the rank of the stoichiometric matrix $N$).  Then, using the notation in the mass-action equation~\eqref{eq:odeYCG}, the network $\mathcal G$ is multistationary if and only if the system
\begin{equation}\label{eq:multi}
Y A_\k x^Y=0,\qquad Z x=c 
\end{equation}
has at least two solutions
$x\in \R^n_{>0}$, 
 for some $\k\in \R^r_{>0}$ and $c\in \R^d$.

There are several methods to evaluate multistationarity quite effectively (see, for instance,~\cite{BanajiCraciun2010, BP-inher, CFMW, PKC, DPST, mss-review, WangSontag} and implementations \cite{crnttoolbox,control}).
One such approach, which precludes multistationarity, is to show that the function defined by the left-hand side of the equations in \eqref{eq:multi} is injective for all parameter values $\kappa \in \mathbb{R}^r_{>0}$ and $c \in \mathbb{R}^d$. If so, we say that the network is \emph{injective}. 

This idea was introduced by Craciun and Feinberg \cite{ME_I} and then further developed by numerous authors, e.g.~in  \cite{BanajiCraciun2009,BanajiCraciun2010,BanajiDonnell,gnacadja_linalg,ME_II,ME_semiopen,MR14,PKC,ShinarFeinberg2012,signs,feliu-bioinfo,FeliuWiuf_MAK,wiuf-feliu-power-law,feliu_newinj,JoshiShiu}. For the purpose of this chapter, we highlight a formulation   given in \cite{FeliuWiuf_MAK} (see also \cite{signs}). 
We begin by removing 
$d$ linearly dependent rows from the stoichiometric matrix $N$, to form a new (full rank) matrix $N'\in \R^{s\times r}$. Let $B\in \R^{n\times r}$ be the matrix of mass-action exponents, that is, the $i$-th column is $y$ if the $i$-th reaction is $y\rightarrow y'$. 
Then, for $\k=(\k_1,\dots,\k_r)$ and $\lambda=(\lambda_1,\dots,\lambda_n)$, we consider the following matrix:
\[ M(\k,\lambda)= \begin{bmatrix}
N' \diag(\k) B^\top \diag(\lambda) \\
Z
\end{bmatrix} \in \R[\k,\lambda]^{n\times n}.\]
With this in place, we have the following result~\cite{FeliuWiuf_MAK}. 

\begin{theorem}[Injectivity criterion] \label{thm:injective}
The network $\mathcal G$ is injective (and hence $\mathcal G$ is non-multistationary) if and only if all coefficients of the polynomial $\det(M(\k,\lambda))$ 
have the same sign, and the polynomial is not the zero polynomial. 
\end{theorem} 

Theorem~\ref{thm:injective} is shown by first establishing that the network is injective if and only if  $\det(M(\k,\lambda))$ never vanishes, and then showing that, unless $\det(M(\k,\lambda))$ is the zero polynomial, this non-vanishing holds if and only if all  coefficients of $\det(M(\k,\lambda))$ have the same sign.  This second equivalence is in turn proven by showing that $\det(M(\k,\lambda))$ is homogeneous in $\k,\lambda$ and is {\em multiaffine}, that is, each exponent is either $0$ or $1$.

\subsection{The human story. }  \label{sec:human-story-injectivity}
Feliu recounts the following:

\emph{Several approaches related to injectivity of networks were unified in \cite{signs}. This happened after several of us -- including Alicia Dickenstein -- realized at the Dagstuhl Seminar on `Symbolic methods for chemical reaction networks' held in November 2012 at Schloss Dagstuhl, Germany, that several of our independent works had the same core techniques and principles.}

\emph{
A few years later, 
while I was visiting Alicia at ICTP in Trieste in March 2015 to work on a common book project, Alicia said something on these lines (which pertains to ideas in the proof of Theorem~\ref{thm:injective} mentioned above): `Bernd told me that the fact that $\det(M(\k,\lambda))$   does not vanish if and only  all its coefficients have the same sign, is simply because the Newton polytope of $\det(M(\k,\lambda))$ is an $n$-dimensional cube and hence all points are vertices.'  
At that time, I was not familiar with the relationship between signs of polynomials and properties of the Newton polytope, although later on Timo de Wolff referred to such a statement as `folklore in real algebraic geometry', and indeed, it goes  back at least to Reznick \cite{reznick}. I was working on another project on multistationarity \cite{CFMW}, which required deciding when a polynomial, which has both positive and negative coefficients, attains a negative value over the positive orthant.  A bit of digging led me to understand this relation (Proposition~\ref{prop:signs}  below), and since then, that result and related ideas are behind many of works in the field \cite{CFMW,torres-feliu,CFM,TelekDescartes,circuits}. }

\subsection{Signs, Newton polytope, multistationarity, and dynamics} 
We now consider multivariate polynomials
\begin{align*}
	p(x)=\sum_{\alpha\in \sigma(p)} a_\alpha x^{\alpha}  ~\in ~ \R[x_1,\dots,x_n],
\end{align*}	
 where 
$\sigma(p)$ is a finite subset of $ \Z_{\geq 0}^n$  and $a_\alpha\neq 0$ (for all $\alpha$). 
The \emph{Newton polytope} of $p$ is the convex hull of the exponent vectors of $p$:
\begin{equation}\label{eq:NP}
\New(p) := {\rm Conv} \big(\{ \alpha\in \sigma(p)\}\big) \subseteq \R^n. 
\end{equation}
Given a face $F$ of the Newton polytope $\New(p)$, 
the polynomial 
$p_F$ is obtained by considering only the exponents in $F$:
\[
p_F(x):= \sum_{\alpha\in F\cap \sigma(p)} a_\alpha x^{\alpha}. \]

With this notation, the result recounted in Section~\ref{sec:human-story-injectivity} says the following:
\begin{proposition}\label{prop:signs}
If there exists  $x^*\in \R^{n}_{>0}$ such that 
$p_F(x^*)<0$, 
then there also exists  $y\in \R^{n}_{>0}$ such that $p(y)<0$.
\end{proposition}

Hence, if $F$ is a vertex $\alpha$, then $p_F(x)=a_\alpha x^\alpha$, and Proposition~\ref{prop:signs}
establishes that if the coefficient of $x^\alpha$ is negative, then $p$ attains  some negative values in the positive orthant. 

The study of the signs that a polynomial attains on the positive orthant leads to answering difficult questions pertaining to reaction networks. 
One of the first results in this direction arose when Craciun, Koeppl, and Pantea analyzed quadratic systems with fixed rate constants $\k$ using the injectivity approach to precluding multistatinonarity described earlier (recall Theorem~\ref{thm:injective})~\cite{PKC}.

Another setting in which Proposition~\ref{prop:signs} can be harnessed 
arises from networks with no boundary steady states  in  stoichiometric compatibility classes that intersect $\R^n_{>0}$, and which satisfy an additional condition.  For systems arising from such networks (and some choice of rate constants $\kappa$), the authors of \cite{CFMW} applied degree theory to 
construct 
a polynomial $p_{\kappa} \in \R[\xi_1,\dots,\xi_q]$ such that\footnote{The polynomial $p_{\kappa}$ refers to the polynomial denoted by $(-1)^s a(\hat x)$ in~\cite{CFMW}.}
\begin{itemize}
\item if $p_{\kappa}(\xi)>0$ for \emph{all}  $\xi \in \R^{q}_{>0}$, then every stoichiometric compatibility class has exactly one positive steady state; 
\item 
if $p_{\kappa}(\xi)<0$  for \emph{some} $\xi \in \R^{q}_{>0}$, 
then there exists a stoichiometric compatibility class with at least two positive steady states. 
\end{itemize}

It follows that the values of $\k$ for which the polynomial $p_{\k}$ attains 
negative values on $\R^q_{>0}$ are precisely those $\k$ that enable multistationarity.  The only subtlety arises when, for some $\kappa$, the polynomial $p_{\kappa}$ takes both positive and zero values; but this situation typically does not arise in applications.
Conditions for multistationarity can therefore be found by quantifier elimination, but this is infeasible in practice.  
Indeed, the polynomial $p_{\k}$ (which we now view as a polynomial in both $\k$ and $\xi$) is typically large. 
Instead, it is often profitable to analyze the coefficients of $p_{\kappa}$, and to apply Proposition~\ref{prop:signs}.

However, Proposition~\ref{prop:signs} 
is silent when all negative coefficients arise from non-vertices of the Newton polytope. 
This situation happens, for instance, with the  distributive phosphorylation mechanism $\mG_1$  in \eqref{eq:G1} considered above.  Feliu recalls how this problem was eventually resolved:  {\em Timo de Wolff and I had talked for a while about using circuit numbers and tools from SONC (Sums Of Non-negative Circuit functions) to determine the parameter region of multistationarity for $\mG_1$. 
Bernd encouraged and brought Timo and me together with his Ph.D. student Nidhi Kaihnsa -- plus 
Timo's Ph.D. student O\u{g}uzhan Y\"ur\"uk -- to address this problem.  Our results are now published in \cite{circuits} and further developed in \cite{circuits2}.
}

\medskip
Another context in which Proposition~\ref{prop:signs} is helpful arises
when studying bistability and Hopf bifurcations, as follows. 
Let $J(x)$ denote 
the Jacobian matrix of $YA_\k x^Y$ in \eqref{eq:odeYCG}. 
For a steady state $x^*$, 
if all eigenvalues of $J(x^*)$ have negative real part, then $x^*$ is {\em asymptotically stable}, and trajectories starting nearby converge to $x^*$ \cite{perko}. If all eigenvalues of $J(x^*)$ have negative real part except for a 
 complex-conjugate pair of purely imaginary eigenvalues -- and an additional non-degeneracy condition holds with respect to moving one parameter value -- then {\em Hopf bifurcations} arise, and as a consequence, periodic solutions exist \cite{perko}.  Periodic solutions can be biologically significant.

Finding eigenvalues of matrices with symbolic entries, like $J(x)$, is not possible in general. However, an alternative approach is to consider the Hurwitz matrix $H$ associated with the characteristic polynomial of $J(x)$, and then to use the Routh-Hurwitz criterion~\cite{Barnett:hurwitz-routh} to characterize the number of 
eigenvalues with positive vs.\ negative real part via the signs of the principal leading minors of $H$. 
 Hence, Proposition~\ref{prop:signs} can again be used to assess stability and Hopf bifurcations; see for instance, ~\cite{Errami,yang-hopf,BifTheo-009,deciding:kahoui}.  
 This approach  has been explored for phosphorylation networks and other signaling networks~\cite{torres-feliu,OSTT,mixed,CFM}.

 \medskip
 Finally, another algebraic approach for finding Hopf bifurcations, based on resultants and the Bezout matrix, was given by Guckenheimer, Myers, and Sturmfels in 1997 \cite{sturmfels:hopf}.

	\section{Counting steady states, mixed volume, and numerical algebraic geometry} \label{sec:counting}

As we saw in the previous section, 
various results and tools are available for assessing whether a given network or family of networks is multistationary. 
Once a network is known to be multistationary, the next problem is to compute -- or at least to bound -- the maximum number of positive steady states (and, ideally, the maximum number of those that are stable).  
This problem of counting (or bounding) the number of positive steady states is an important and active research topic.

In this section, we highlight two of Sturmfels's contributions in this direction.  Specifically, his work used for the first time (in this area) two concepts from adjacent areas -- mixed volume (Section~\ref{sec:MV}) and numerical algebraic geometry (Section~\ref{sec:NAG}).  Both approaches are now quite active areas of research.  
At the end of the section, we highlight additional methods for counting steady states (Section~\ref{sec:more-approaches-count}).

\subsection{Mixed volume} \label{sec:MV}
One way to bound the number of positive steady states, that is, the number of positive roots of the corresponding parameterized polynomial system, is to determine the maximum number of roots in the complex torus $\C^*$.  One is then led to the concept of {\em mixed volume} from convex geometry~\cite{betke}.  We do not define this precisely here, but we note that the mixed volume of polytopes is readily computable using computer algebra software, and its key property is 
the content of Bernshtein's theorem~\cite{bernstein}, as follows. 

\begin{proposition} \label{prop:bernstein}
Let $g_1,  \dots, g_n \in \mathbb{R}[x_1,  \dots, x_n] $.
Then, counted with multiplicity, 
the number of isolated solutions in $(\mathbb{C}^*)^n$ of the system 
$g_1(x) =   \cdots = g_n(x) = 0$ is at most the mixed volume of the Newton polytopes $\New(g_1), \ldots , \New(g_n)$.
\end{proposition}

Therefore,
an upper bound on the number of positive steady states can be obtained from 
the mixed volume of an appropriate set of Newton polytopes.  
Several ways to build such polytopes from the right-hand sides of the ODEs are described below.

In~\cite{case-study}, Sturmfels and collaborators were the first to put related ideas into action. 
However, it is important to note that {\em the mixed volumes they define were \emph{not} constructed to be upper bounds on the number of positive steady states}; their motivation was instead to assess which of the many parametrizations arising from bases of an algebraic matroid, as mentioned in Section~\ref{sec:steady-state-params}, are the `best'\footnote{The idea is that parametrizations giving rise to mixed volumes that are equal to the true number of solutions in  $(\mathbb{C}^*)^n$ generate good start systems for applying numerical homotopy-continuation methods.}.  
Nevertheless, these mixed volumes often appear to be good (but not tight) upper bounds.  
For instance, these mixed volumes range from $5$ to $45$ for a certain {\em Wnt shuttle network} (a model for an important biological signaling pathway, with $19$ species and $36$ reactions), as shown in Table~\ref{tab:mv}, while the  
conjectured maximum number of positive steady states is $3$~\cite{case-study}.

\begin{table}[ht]
\caption{Three reaction networks, their maximum number of positive steady states, matroidal mixed volume, augmented-system mixed volume, and references for these results.  Both `$\geq 3$' entries are conjectured to be `3'~\cite{case-study, OSTT}.  The entry `$5-45$' refers to the fact that the mixed volumes obtained from various choices of parametrizations range from $5$ to $45$.  The `?'' entries were not computed, due to the complicated nature of matMV that makes it difficult to apply to general networks. }\label{tab:mv}
\begin{center}
\begin{tabular}{l c c c c}
\hline
	Network & Max \# & matMV & augMV & Reference\\
			& positive  & & \\
			&  steady states & & \\
	\hline
Wnt & $\geq 3$ & $5-45$ & 56 & \cite{case-study, Obatake-thesis} \\
ERK & $\geq 3$ & ? & 7 & \cite{OSTT} \\
Fully irreversible ERK & 1 & ? & 3 & \cite{OSTT}\\
\hline
\end{tabular}
\end{center}

\end{table}%

With the door opened to mixed-volume computations, other formulations of mixed volumes for networks (that is, using other associated polytopes) were proposed~\cite{gross-hill, OSTT, Obatake-thesis, mv-small-networks}.  
These mixed volumes are compared in the recent dissertation of Obatake~\cite{Obatake-thesis}, some of which we summarize here.  Sturmfels' original mixed volume, which Obatake calls the {\em matroidal mixed volume} (matMV), is somewhat complicated.  Indeed, it involves minimizing over choices of generators of a saturation ideal coming from the basis of an algebraic matroid, and this minimum depends on the choice of the basis (see Table~\ref{tab:mv}). 

In contrast, the mixed volume of Obatake, Shiu, Tang, and Torres~\cite{OSTT} -- the {\em augmented-system mixed volume} (augMV) -- is straightforward to define: it is simply the mixed volume of the Newton polytopes of the polynomials obtained from the right-hand sides of the mass-action ODEs after linearly redundant ODEs are replaced using conservation laws. The mixed volume of Gross and Hill~\cite{gross-hill} is similar in spirit, and even lends itself to combinatorial arguments, thereby enabling mixed-volume computations for infinite families of networks~\cite{gross-hill}.  Moreover, both augMV and the Gross--Hill mixed volume fulfill the inequality, mentioned earlier, that motivates our interest in the mixed volume~\cite{Obatake-thesis}:
\begin{align} \label{eq:ineq-MV-1}
\mathrm{Max}~\#~\mathrm{positive~steady~states} \quad \leq \quad 
\mathrm{ augMV} \quad \leq \quad
 \text{ Gross--Hill~MV}
\end{align}

In Table~\ref{tab:mv}, 
we compare the matroidal and augmented-system mixed volumes for the Wnt shuttle network.
Also shown in that table is data on the 
 {\em ERK network} (and also an irreversible version) from~\cite{long-term,OSTT}.  For the ERK network, 
the mixed-volume bound is $7$, while the maximum number of positive steady states is conjectured to be $3$~\cite{OSTT}.

We now consider a fourth type of mixed volume, which Obatake calls the steady-state-parametrization mixed volume (sspMV)~\cite{Obatake-thesis}. 
This mixed volume, which can be viewed as a saturation-free version of Sturmfels's matMV, is defined for networks having a steady-state parametrization.  Also, like augMV and the Gross--Hill mixed volume, it bounds the number of positive steady states:
\begin{align} \label{eq:ineq-MV-2}
\mathrm{Max}~\#~\mathrm{positive~steady~states} \quad \leq \quad 
\mathrm{ sspMV} \ ~.
\end{align}
The procedure for computing sspMV (given some parametrization) is to substitute the parametrization into the conservation laws, clear denominators, and then compute the mixed volume of the resulting Newton polytopes.  We illustrate this procedure in the following example.

\begin{example}
We return to the McKeithan network \eqref{eq:mckeithan} and the positive steady-state parametrization in \eqref{eq:mckeithanpar}. 
We plug this parametrization into the following conservation laws, as in \eqref{eq:mckeithan_cons}: $x_1+x_3+x_4 - c_1 =0$ and $x_2+x_3+x_4- c_2=0$, which yields the following equations in 
$x_1,x_2$:
\[
x_1+\tfrac{\k_1}{\k_2+\k_3} x_1x_2+\tfrac{\k_1\k_3}{(\k_2+\k_3)\k_4} x_1x_2 - c_1=0~,\qquad x_2+\tfrac{\k_1}{\k_2+\k_3} x_1x_2+\tfrac{\k_1\k_3}{(\k_2+\k_3)\k_4} x_1x_2-c_2=0~. \]
For these two polynomials, the Newton polytopes are (respectively) the following triangles, which we denote by $P_1$ and $P_2$:

\begin{center}
\begin{tikzpicture}[scale=1.2] 
\node[left] (a) at (0,0) {$(0,0)$ };
\node[right] (a) at (1,0) {$(1,0)$ };
\node[right] (a) at (1,1) {$(1,1)$ };
\draw[fill=black!10!white,line width=0.8pt]  (0,0) -- (1,0) -- (1,1) -- (0,0);
\node[left] (a) at (1,0.3) {$P_1$ };
\draw[fill=black] (0,0) circle (2pt);
\draw[fill=black] (1,0) circle (2pt);
\draw[fill=black] (1,1) circle (2pt);
%
\node[left] (a) at (4,0) {$(0,0)$ };
\node[left] (a) at (4,1) {$(0,1)$ };
\node[right] (a) at (5,1) {$(1,1)$ };
\draw[fill=black!10!white,line width=0.8pt]  (4,0) -- (4,1) -- (5,1) -- (4,0);
\node[left] (a) at (4.6,0.7) {$P_2$ };
\draw[fill=black] (4,0) circle (2pt);
\draw[fill=black] (4,1) circle (2pt);
\draw[fill=black] (5,1) circle (2pt);
\end{tikzpicture}
\end{center}
The mixed volume $MV(P_1,P_2)$ equals $2$, and so 
${\rm sspMV}=2$ (for the parametrization we used). 
For this model, the mixed-volume bound~\eqref{eq:ineq-MV-2}
 is not tight (recall that there is always a unique positive steady state -- and it is complex-balanced -- as, from~\eqref{eq:mck-def}, the deficiency is 0).  
\end{example}

In light of the bounds~\eqref{eq:ineq-MV-1} and~\eqref{eq:ineq-MV-2}, it is natural to compare augMV and sspMV.  In fact, 
we can have 
$\text{augMV} \leq \text{ sspMV}$ 
for some steady-state parametrization, while 
$\text{augMV} > \text{ sspMV}$ 
for some `better' ones (and so sspMV improves the bound)~\cite{Obatake-thesis}.  
Our ideas have therefore come full circle: finding the `best' parametrizations was the original motivation for Sturmfels to propose mixed-volume computations for reaction networks!

\subsection{Numerical algebraic geometry} \label{sec:NAG}
As noted above, the Wnt shuttle network is conjectured to admit at most three positive steady states.  In the first application of {\em numerical algebraic geometry} to reaction networks, Sturmfels and collaborators gave numerical evidence for this conjecture~\cite{case-study}.  Specifically, they used the software {\tt Bertini}~\cite{BHSW06} to sample 10,000 pairs of parameter and total-amount values, and checked that the number of positive steady states was never more than 3.  
They also performed a similar analysis to explore how robust multistationarity of this network is to perturbations in the parameter and total-amount values.

Numerical algebraic geometry continues to play an important role in investigations of biochemical reaction networks
\cite{NAG-life-sci, bates-gunawardena, harrington-decompose}.  For details on what  numerical algebraic geometry is all about, we refer the reader to another chapter in this volume, `Numerical Nonlinear Algebra'~\cite{bates2023numerical}.

 \subsection{Other approaches to counting steady states} \label{sec:more-approaches-count}
Here we describe additional approaches to counting or bounding the number of positive steady states. 
For specific families of biochemical reaction networks, several authors established unlimited multistability  , that is, these families have no upper bound on the number of stable, positive steady states \cite{kothamachu-unlimited,ThomGuna09a,rendall-unlimited}.
Next, the methods of Dickenstein, in collaboration with 
Fr\'ed\'eric Bihan and 
 Magal\'i Giaroli, 
 rely on finding `positively decorated' simplices from within the convex hull of all monomials appearing in the mass-action ODEs~\cite{bihan-dickenstein-giaroli-lower,Giaroli-Bihan-Dickenstein}.   
Finally, Nida Obatake and Elise Walker go beyond the mixed volume to volumes of `Newton-Oukonkov bodies' to establish improved bounds on the number of steady states~\cite{NO-bodies-obatake-walker}.

	\section{Convex hulls of trajectories and the attainable region problem} \label{sec:attainable}

\subsection{The human story} Shiu recalls the following:

 {\em 
The story begins in the early 1990s, when Martin Feinberg reached out to Victor Klee, who was one of Bernd's PhD advisors, for some help proving a result concerning boundary points of convex sets.  Years later, in 2009, Martin mentioned these interactions to some of us, including Bernd, at a conference, and subsequently sent some notes and references on what is called the 
`attainable region problem'.  

This problem is of significant industrial interest and also inherently involves optimization and convex geometry.  Bernd therefore hoped to generate interest in this problem among the convex/real algebraic geometers at Berkeley at the time, such as graduate student Cynthia Vinzant, so he had me give an informal talk on the topic.  

Nothing came of it until years later.  
As Bernd wrote in an email to one of his PhD students, Nidhi Kaihnsa:
 ``I had a very exciting discussion with a famous Chemical Engineer [Larry Biegler] yesterday, along with Cynthia Vinzant. The outcome is a very nice research theme I would like to suggest to you as [a] possible dissertation topic. The scientific context is `Attainable Region Theory'. Please google this. It looks very non-mathematical, right? But, I am absolutely convinced that this is a `goldmine' for someone like you, who is open to the natural sciences. It's all about Convex Algebraic Geometry.'' Of course, Bernd's assessment proved spot-on, and the results Nidhi and Bernd achieved are described below.
}

\subsection{The attainable region problem}
 {\em Attainable region theory} is an important research area in chemical engineering; it pertains to the mathematics of understanding
two related sets, the {\em attainable region} and the convex hull of trajectories. 
Here, trajectories $\{\phi(t) \mid t \geq 0\}$ are images of solutions $\phi(t)$ to~\eqref{eq:ode}.  The convex hull of a trajectory is the smallest convex set containing the trajectory, and the {\em attainable region of a trajectory} (or the initial condition of the trajectory) is the (generally larger) set that is the smallest convex set that not only contains the trajectory but also is forward-closed with respect to the dynamics.  

More generally, an {\em attainable region} is a convex set that describes all possible concentration vectors that can be obtained by basic chemical processes (mixing, reacting, and inflows) for a given  reaction network and inflow rates.  The goal, in chemical engineering applications, is to characterize (the boundary of) the attainable region, in order to optimize over this region.  For instance, one might wish to maximize the production of a desired chemical~\cite{ming2016attainable}. 
As recounted by Feinberg~\cite{feinberg2002toward}, 
the attainable region problem was initiated by Horn in 1964~\cite{horn1964attainable}. Significant results were elucidated by Feinberg, Crowe, Glasser, and Hildebrandt from the 1980s onward~\cite{FeinH_ORD1, hildebrandt1990geometry}.

\subsection{Contributions of Sturmfels and collaborators}
The main results of Sturmfels and co-authors are as follows.  First, Kaihnsa
proved that for chemical reaction networks that are `linear', that is, the ODEs are given by linear polynomials, the two regions of interest, the attainable region and the convex hull of a trajectory, are always equal~\cite{kaihnsa}.  Moreover, under some hypotheses, this region is a `spectrahedral shadow' (a convex set that can be described by linear matrix inequalities).  Finally, Kaihnsa conducted some numerical experiments and conjectured that for a certain class of networks, the attainable region and the convex hull of a trajectory are again equal.

This conjecture was resolved (in the negative) in the follow-up paper of Kaihnsa and Sturmfels, with Ciripoi and L\"ohne~\cite{CKLS}.  
The counterexample, and many other related examples, were enabled by an algorithm that the authors developed and implemented.  Specifically, the algorithm focuses on decomposing into `patches' the boundary of the convex hull of a trajectory. 
Moreover, the mathematical foundation underlying these results is motivated by an important problem concerning convex bodies, which also is of general interest in the area of convex algebraic geometry: approximating convex hulls by convex polytopes.  Indeed, in this context, the idea of patches was further developed by Plaumann, Sinn, and Wesner~\cite{PSW}.

\section{Closing stories} \label{sec:end}

In this chapter we did not manage to cover all works in the field of reaction networks that arose from matchmakings engineered by Sturmfels. 
For instance, 
recent works of 
Brustenga, Craciun, and Sorea~\cite{CM,BCM} 
began while all were at 
MPI Leipzig, where Sturmfels is a director, and Sturmfels is acknowledged for bringing the team together. 
Another example involves Martin Helmer, who was Sturmfels's postdoc when he attended the global attractor conjecture workshop mentioned in Section~\ref{sec:GAC}, and subsequently collaborated with Michael Adamer on model testing and multistationarity for systems with toric steady states~\cite{adamer-helmer-model-testing,adamer-helmer}
before joining Feliu's research team~\cite{feliu-helmer}.

\medskip

Feliu and Shiu share some final thoughts:

\textit{  
Both authors of this chapter have benefited from the support and advice of Bernd, and, had it not been for this, our research careers would look quite different. We are extremely thankful for that. We thought we would honor Bernd by writing the chapter in his style (more or less, as we do need to sleep...). We worked intensively on the first draft for a week, taking advantage of the time zone difference such that the chapter never slept! 
And we must admit that it worked quite well. 
}

\textit{
We hope that with this chapter we managed to convey our appreciation and admiration of 
Bernd and our recognition of his considerable influence on our research field and on all of us in this area.  
We definitely enjoyed writing it and digging into the past. 
}

\begin{center}
{\bf Happy birthday, Bernd!}
\end{center}

\subsection*{Acknowledgements} {\small
The authors acknowledge the IMA 
for permission to use the photo in Figure~\ref{fig:IMA}.  
The authors thank 
Tim Bates, 
Aatmun Baxi, 
Cashous Bortner, 
Crystal Farris, 
Seth Gerberding, 
Nidhi Kaihnsa,
Tung Nguyen, 
Nida Obatake, 
Alexander Ruys de Perez,
and M\'at\'e Telek, 
for helpful comments on an earlier version of this chapter.  
The authors also acknowledge Elizabeth Gross for useful discussions, and several people who shared personal recollections with us: Alicia Dickenstein, Thomas Kahle, and Nidhi Kaihnsa.
Finally, the authors are grateful to three anonymous reviewers, whose detailed comments improved this work.
}

\bibliographystyle{amsalpha}
\newcommand{\etalchar}[1]{$^{#1}$}
\providecommand{\bysame}{\leavevmode\hbox to3em{\hrulefill}\thinspace}
\providecommand{\MR}{\relax\ifhmode\unskip\space\fi MR }
\providecommand{\MRhref}[2]{%
  \href{http://www.ams.org/mathscinet-getitem?mr=#1}{#2}
}
\providecommand{\href}[2]{#2}

\end{document}